\documentclass{article}
\usepackage{geometry,amsmath,amssymb,enumerate,cite,graphicx,caption,url,xcolor,mathrsfs,subfigure}
\usepackage[ruled]{algorithm2e}
\usepackage[thmmarks,amsmath,thref]{ntheorem}
\usepackage{hyperref}
\newtheorem{theorem}{Theorem}[section]
\newtheorem{definition}{Definition}[section]
\newtheorem{lemma}{Lemma}[section]
\newtheorem{assumption}{Assumption}[section]

\newtheorem*{remark}{Remark}
\theorembodyfont{\normalfont}

\theoremsymbol{{$\square$}}
\newtheorem*{proof}{Proof}

\providecommand{\keywords}[1]
{
  \small
  \textbf{Keywords:} #1
}

\providecommand{\MSC}[1]
{
  \small
  \textbf{Mathematics Subject Classification 2020:} #1
}

\numberwithin{equation}{section}
\title{Fictitious Play for Mean Field Games with Optimal Stopping: Convergence and Computation}
\author{Yifan Luo\footnotemark[1], \and Chengfeng Shen\footnotemark[2], \and Jiajun Tong\footnotemark[3], \and Zhennan Zhou\footnotemark[4]}

\begin{document}
\maketitle

\renewcommand{\thefootnote}{\fnsymbol{footnote}}
\footnotetext[1]{School of Mathematical Sciences, Peking University, Beijing, 100871, China. (luoyf@pku.edu.cn)}
\footnotetext[2]{School of Mathematical Sciences, Peking University, Beijing, 100871, China. (Shencf1999@pku.edu.cn)}
\footnotetext[3]{Beijing International Center for Mathematical Research (BICMR), Peking University, Beijing, 100871, China. (tongj@bicmr.pku.edu.cn)}
\footnotetext[4]{Institute for Theoretical Sciences, Westlake University, Hangzhou, 310030, China. (zhouzhennan@westlake.edu.cn)}

\begin{abstract}

 This paper studies mean field games with optimal stopping time (OSMFGs) where agents make optimal exit decisions. Such a model features a coupled obstacle problem and Fokker-Planck equation posing challenges on top of classic mean field games. The nonconvex nature of exit decisions renders the existence of a classic pure strategy equilibrium infeasible, necessitating the consideration of more complex mixed strategy equilibria. 
 
 This paper proposes a generalized fictitious play algorithm that computes OSMFG mixed equilibria by iteratively solving pure strategy systems, i.e., approximating mixed strategies through averaging pure strategies according to a certain updating rule. The generalized fictitious play allows for a broad family of learning rates and the convergence to the mixed strategy equilibrium can be rigorously justified. The algorithm also incorporates efficient finite difference schemes of the pure strategy system. Numerical experiments demonstrate the effectiveness of the proposed method in robustly and efficiently computing mixed equilibria for OSMFGs.
\end{abstract}

\MSC{91A13, 60G40, 65M06.}

\keywords{mean field games, optimal stopping, obstacle problem, fictitious play, finite difference method.}

\section{Introduction}
Recently, large-scale strategic interactions involving many agents, such as online auctions and online voting platforms, have become more prevalent compared to small-scale game scenarios with few players.
Introduced in the seminal works by Lasry and Lions \cite{MFG1,MFG2,MFG3}, the theory of mean field games (MFGs) provides a framework for modeling games with a large number of players. In an MFG model, a continuum of identical players controls the velocities of their own trajectories based on anticipating the evolution of the population distribution, with the goal of maximizing their expected payoffs. The key challenge is to establish the Nash equilibrium, where each agent's anticipated population dynamic aligns exactly with the actual population dynamic resulting from their optimal controls. 
In the Nash equilibrium, each agent's value function satisfies a backward Hamilton-Jacobi-Bellman (HJB) equation, while the population distribution satisfies a forward Fokker-Planck equation. The resulting mathematical structure is a coupled system of forward-backward partial differential equations (PDEs), which provides a complete characterization of the equilibrium.

Beyond the classical MFG settings, recent theoretical advances have extended the MFG framework to incorporate optimal stopping (OSMFG) \cite{osmfgpde}. Instead of controlling the drift term of a stochastic differential equation (SDE) as in a standard MFG, agents in an OSMFG decide whether to exit the game at each state-time point to maximize expected utility. This framework naturally models economic scenarios like firm exit decisions \cite{produce}, or industry turnover \cite{electricity}. Mathematically, the introduction of optimal stopping control transforms the PDE satisfied by the value function from an HJB equation to an obstacle equation.  

To understand this change, consider an individual agent whose state evolves according to a stochastic differential equation with generator $\mathcal{L}$. The agent's goal is to choose a stopping time $\tau$ to minimize an expected cost functional, which consists of a running cost $f(m)$ accumulated over time, and a terminal exit cost $\psi(m)$ incurred upon leaving the game. Both costs depend on the population distribution $m$, capturing the mean-field interaction. Analogous to the classical MFGs, the system of coupled forward‑backward PDEs characterising the pure strategy Nash equilibrium for OSMFGs can be written as:
\begin{equation}\label{PDE4OSMFG_intro}
    \begin{cases}\
\max(-\partial _tu-\mathcal{L} u-f(m),u-\psi(m))=0,&(x,t)\in \Omega\times(0,T);\\
\partial_tm-\mathcal{L}^* m=0,&(x,t)\in\{u<\psi\};\\
m=0,&(x,t)\in\{u=\psi\};\\
u(x,T)=\psi(x,T,m),\,m(x,0)=m^0(x)&x\in \Omega;\\
\text{proper B.C.}
\end{cases}.
\end{equation}
Here, $u$ denotes the value function of an individual player, and $m^0$ denotes the initial distribution of states at time $0$. However, unlike the classical mean field games where agents exert continuous control over their drift velocity, the control in an OSMFG is strictly binary: either stay in the game or exit. This singular control structure means the optimization landscape lacks the necessary convexity to guarantee the existence of a pure strategy Nash equilibrium. A more complex mixed‑strategy Nash equilibrium becomes necessary \cite{osmfgpde}.
In a mixed‑strategy equilibrium, an agent chooses a probability (between 0 and 1) of leaving or staying, rather than a deterministic decision. The corresponding PDE system is
\begin{align}\label{PDE4mix_intro}
\begin{cases}
\max(-\partial_tu-\mathcal{L} u-f(m),u-\psi(m))=0, &(x,t)\in \Omega\times(0,T);\\
\partial_tm-\mathcal{L}^* m=0, &(x,t)\in \{u<\psi\};\\
\partial_tm-\mathcal{L}^* m\le0,\, m\ge 0, &(x,t)\in\Omega\times(0,T);\\
\int_{\{u=\psi\}}(f(m)+(\partial_t+\mathcal{L})\psi(m))m\,dxdt=0.\\
u(x,T)=\psi(x,T,m),\,m(x,0)=m^0(x)&x\in \Omega;\\
\text{proper B.C.}
\end{cases}.
\end{align}
Compared with the pure strategy system \eqref{PDE4OSMFG_intro}, the condition $m=0$ in the region $\{u=\psi\}\times(0,T)$ is relaxed to the conditions $\partial_tm-\mathcal{L}^* m\le0$ and $m\ge 0$, along with the complementary condition:
$$
\int_{\{u=\psi\}}(f(m)+(\partial_t+\mathcal{L})\psi(m))m\,dxdt=0.
$$
The existence of a solution to \eqref{PDE4mix_intro} was first established in \cite{osmfgpde}. However, as indicated in \cite{measureflow}, the proof in \cite{osmfgpde} is flawed\footnote{The weak convergence of $m^{\epsilon}$ in $L^2((0,T),H_0^1(\Omega))$ is not sufficient to conclude the convergence of $\int f(m^\epsilon)m^\epsilon$ when $f$ is not linear.}. Hence, a different approach based on measure flow, rather than PDEs, was developed in \cite{measureflow} to establish the existence of mixed-strategy Nash equilibria. In general, OSMFGs pose greater analytical and computational challenges than the classical MFGs due to the singular nature of the stopping control.

A variety of numerical methods has been proposed for solving MFG models, with early work by Achdou and Capuzzo-Dolcetta introducing finite difference techniques and analyzing their convergence properties \cite{achdou_numericalmfg,achdou_numericalmfg2,achdou_numericalmfg_conv}, and since then numerous other approaches have been developed—a comprehensive summary of which is provided by Laurière \cite{numericalMFG}. However, few algorithms have been proposed so far for OSMFG. There have been two main types of numerical algorithms developed for this problem. One is Uzawa's algorithm, proposed in \cite{Uzawa}, which operates under the assumption of monotone running costs. 
This optimization-based method handles the mixed strategy equilibrium PDEs system raised in \cite{osmfgpde}, but is limited to the stationary case. However, in many practical applications, there is a need to solve time-dependent problems, which cannot be directly tackled by the Uzawa method. 
The other method utilizes a large-scale linear programming approach with fictitious play, as introduced in \cite{LP}. This algorithm relies on the measure flow formulation of OSMFG presented in \cite{measureflow} and has been applied to practical problems such as modeling games in electricity markets \cite{electricity}. However, a major drawback of this approach is the need to solve a large-scale linear program at each iteration, which reduces computational efficiency and robustness. Additionally, while this method can provide information on the crowd propagation, it does not reveal the individual stopping strategies employed by each player. In summary, efficient numerical methods for solving general OSMFG problems remain underdeveloped due to the complex mathematical structure of such games. Further algorithm development is needed to handle the obstacle equations and mixed strategy equilibria inherent to OSMFGs in a computationally practical manner.

In this paper, we propose an iterative algorithm for computing the solution to the mixed strategy equilibrium PDEs system \eqref{PDE4mix_intro} for potential OSMFGs, by using a generalized fictitious play to the pure strategy equilibrium PDEs system \eqref{PDE4OSMFG_intro}. Fictitious play is a fixed-point iteration with a learning rate $\delta_n=1/n$, which has been applied to find equilibria for potential classical MFG models \cite{fictitiousplay}. In this paper, we generalize the notion of fictitious play by requiring the learning rate $\delta_n$ to satisfy the Robbins-Monro condition
\begin{align*}
  \textstyle  \sum_n\delta_n=\infty \quad \mbox{and} \quad \sum_n\delta_n^2<\infty.
\end{align*} 
The core idea of this algorithm is to repeatedly solve pure strategy equilibrium systems and construct the mixed strategy via convex combinations of pure strategies, thereby avoiding the complexity of directly handling the mixed strategy equilibrium system. In each iteration, given the current average distribution $\bar{m}_n$, we first solve the backward problem
\begin{equation}\label{eq:u_fictitious_intro}
\begin{cases}
\max(-\partial_t u_{n+1} - \Delta u_{n+1} - f(\bar{m}_n),\; u_{n+1} - \psi(\bar{m}_n)) = 0, & (x,t)\in \Omega\times(0,T);\\
u_{n+1}(x,T) = \psi(x,T,\bar{m}_n), & x\in\Omega;\\
\text{proper B.C.}, & x\in\partial\Omega
\end{cases}
\end{equation}
to obtain a new value function $u_{n+1}$. We then solve the forward problem
\begin{equation}\label{eq:m_fictitious_intro}
\begin{cases}
\partial_t m_{n+1} - \Delta m_{n+1} = 0, & (x,t)\in\{u_{n+1}<\psi\};\\
m_{n+1} = 0, & (x,t)\in\{u_{n+1}=\psi\};\\
m_{n+1}(x,0) = m^0(x), & x\in\Omega;\\
\text{proper B.C.}, & x\in\partial\Omega
\end{cases}
\end{equation}
to obtain a new proposed distribution $m_{n+1}$. The average distribution is then updated via a mixture of $\bar{m}_n$ and $m_{n+1}$ with learning rate $\delta_{n+1}$:
\[
\bar{m}_{n+1} = \delta_{n+1} m_{n+1} + (1-\delta_{n+1})\bar{m}_n.
\]
When $\delta_n = 1/n$, $\bar{m}_n$ becomes the running average of the proposed distributions, where each proposed distribution is a pure strategy from the corresponding iteration. Hence $\bar{m}_n$ can be interpreted as a mixed strategy. In practice, both the backward and forward problems can be solved with standard numerical PDE tools (e.g., finite differences), making the algorithm easy to implement.

\paragraph{Contributions.}
The main contributions of this paper are twofold.

\emph{Theoretical contributions.} We prove the strong compactness in $L^2((0,T),L^2(\Omega))$ of the sequence $\{\bar{m}_n\}$ generated by the generalized fictitious play, based on the decreasing property $\partial_t \bar{m}_n - \Delta \bar{m}_n \le 0$. Using this compactness, we establish convergence of a subsequence of $(u_n,\bar{m}_n)$ to a mixed‑strategy Nash equilibrium of \eqref{PDE4mix_intro} for potential games\footnote{Without the potential game structure, the same approach applied to the sequence $m^{\epsilon}$ in \cite{osmfgpde} also yields strong compactness, thereby fixing the existence proof in \cite{osmfgpde}.}. This result theoretically bridges the pure and mixed strategies, enabling the computation of mixed equilibria by solving only the pure‑strategy systems at each iteration. Furthermore, unlike the measure-flow based optimization methods which only track macroscopic crowd propagation, our formulation inherently constructs the exact probabilistic individual stopping strategies by averaging the sequence of pure strategies.
Similar idea first appeared in Section~6 of \cite{measureflow}, but that earlier work assumed a very restricted form of the running cost $f$ due to lack of strong compactness, and the learning rate was chosen as the minimizer of the potential along the update direction—a choice less flexible than our generalized learning rate. Additionally, we analyze the fully discretized system. Various numerical methods exist for obstacle equations, including $\mathbb{L}^1$ optimization \cite{obstacle_L1_method}, finite element schemes \cite{obstacle_FEmethod}, and several finite difference schemes \cite{obstacle_PLS, obstacle_MultigridFD, obstacle_SOR}. We adopt a finite difference method that accurately captures the free boundary while remaining simple to implement. For the fully implicit finite difference scheme, we prove that the sequence $(u_n,\bar{m}_n)$ generated by the generalized fictitious play converges to the fully discretized mixed‑strategy equilibrium system.

\emph{Numerical contributions.} We develop an efficient semi‑implicit finite difference scheme to solve the obstacle equation and the Fokker–Planck equation, and demonstrate its effectiveness through two numerical experiments. The first experiment considers a one-dimensional non-local OSMFG that admits a pure strategy equilibrium, and serves to illustrate the convergence behavior of the algorithm as well as the flexibility of the generalized learning rate in this relatively simple setting. The second experiment considers a two-dimensional local OSMFG for which no pure strategy equilibrium exists, and demonstrates that the condition on $\delta_n$ stated above is indispensable for convergence to the mixed strategy equilibrium in such cases.

The rest of this paper is organized as follows. In Section 2, we introduce the OSMFG model and present the PDE systems characterizing the pure and mixed strategy equilibria. In Section 3, we prove the convergence of generalized fictitious play in the continuous setting and construct the fully-implicit and semi-implicit schemes for the obstacle problem and the Fokker-Planck equation. We also prove the convergence of generalized fictitious play combined with the fully-implicit obstacle scheme. In Section 4, we present numerical experiments validating our algorithms. Finally, in Section 5, we summarize the paper's main results and discuss potential extensions.

\section{Model}
In this section, we introduce the PDE system characterizing the mixed strategy equilibrium in mean field games with optimal stopping (OSMFG) and establish a compactness result that is crucial for the existence of solutions to this system. We begin by presenting the single agent optimal stopping problem, along with the associated PDE for the value function. We then proceed to develop the OSMFG framework. While a pure strategy equilibrium can be defined analogously to that in the classical mean field games, such an equilibrium may fail to exist due to the singular stopping controls. To address this, we relax the constraint on agents' choices and define a mixed strategy equilibrium as in \cite{osmfgpde}. Finally, we review the existence result of this mixed strategy equilibrium in \cite{osmfgpde}\footnote{As stressed in the introduction, the proof in \cite{osmfgpde} is flawed.}, and introduce a compactness argument that completes the proof.

\subsection{Mean Field Game with Optimal Stopping}
\subsubsection{Single Agent Optimal Stopping and Obstacle Equation}
Let $X_t$ be the state variable of an agent that satisfies a stochastic differential equation (SDE): $$dX_t=b(X_t,t)dt+\sigma(X_t,t) dW_t,$$ where $W_t$ is a standard Brownian Motion. The only admissible control is choosing the time to exit the game. Let $f:\mathbb{R}^n\times \mathbb{R}^+\rightarrow\mathbb{R}$ be the running cost function and $\psi:\mathbb{R}^n\times\mathbb{R}^+\rightarrow \mathbb{R}$ be the stopping cost. Given $T>0$, the optimal stopping problem can be formulated as
$$
\min_{0\le \tau \le T}\mathbb{E}\left[\int_{0}^{\tau}f(X_s,s)ds+\psi(X_\tau,\tau)|X_0=x\right],
$$
where $\tau$ is the stopping time which is adapted to the filtration generated by $W_t$. Define 
$$u(x,t)=\min_{t\le\tau\le T}\mathbb{E}\left[\int_{t}^{\tau}f(X_s,s)ds+\psi(X_\tau,\tau)|X_t=x\right]$$
as the value function. Assuming that $u$ is smooth enough, we can derive the PDE system satisfied by $u$ formally by the dynamic programming principle. Details of derivation are omitted here.  Let $A(x,t)=\sigma(x,t)\sigma(x,t)^T=(a_{ij}(x,t))$ be the diffusion matrix, and $\mathcal{L}:=b(x,t)\cdot \nabla+\frac{1}{2}\sum_{i,j}a_{ij}(x,t)\partial_{ij}$ be the generator of the diffusion process defined by the SDE satisfied by $X_t$. The governing PDE system for $u$ is as follows:
\begin{equation}\label{obstacle}
\begin{cases}
\max(-\partial _tu-\mathcal{L} u-f,u-\psi)=0,&(x,t)\in \mathbb{R}^n\times(0,T);\\
u(x,T)=\psi(x,T),&x\in \mathbb{R}^n.\\
\end{cases}
\end{equation}

 An optimal control for an agent is to leave the game at $$\tau=\min_{0\le t\le T}\{u(X_t,t)=\psi(X_t,t)\}.$$ That means an agent would definitely  leave the game at time $t$ and state $x$ if and only if $(x,t)\in\{(x,t):u(x,t)=\psi(x,t)\}$ (or $\{u=\psi\}$ for simplicity). We note that this might not be--- the unique optimal control for the optimal stopping problem: in some regions of $\{u=\psi\}$, agents may be indifferent between exiting and staying, since both actions yield the same cost. This indifference naturally leads to the mixed-strategies considered later. Problem~(\ref{obstacle}) is a parabolic obstacle problem---a type of free-boundary problem in which the solution $u$ satisfies a PDE on a domain whose boundary is determined by the solution itself. Obstacle problems have been extensively studied since the 1980s; we refer the reader to Chapter~5 of \cite{obstacle} for further details.
\subsubsection{Mean Field Games with Optimal Stopping} 
In mean field games with optimal stopping (OSMFGs), a continuum of identical agents make optimal stopping decisions to minimize their expected costs, with each agent's cost depending on the collective behavior of the crowd. The goal is to characterize the symmetric Nash equilibrium, where all agents adopt the same policy and each agent's policy is the optimal response to the crowd distribution it induces.

Such an equilibrium can be constructed via a fixed-point consistency condition. Agents first postulate a crowd distribution over time; based on this guess, each agent solves an optimal stopping problem, yielding an actual crowd distribution. An equilibrium is reached when the actual distribution coincides with the guessed one.

Let $T>0$ and $\Omega\in\mathbb{R}^n$ be a bounded smooth region. Let $$dX_t=b(X_t,t)dt+\sigma(X_t,t) dW_t$$ be the dynamic of a representative agent. Let $\mathcal{M}(\Omega\times[0,T])$ denote the space of finite Borel measures on $\Omega\times[0,T]$.\footnote{Later we will work with the subset of  measures with $L^2((0,T);H_0^1(\Omega))$ densities.}  Let $f:\Omega\times[0,T]\times\mathcal{M}(\Omega\times[0,T])\rightarrow\mathbb{R}$ be the running cost function and $\psi:\Omega\times[0,T]\times\mathcal{M}(\Omega\times[0,T])\rightarrow\mathbb{R}$ be the stopping cost. Given the crowd distribution $\mu \in \mathcal{M}(\Omega\times[0,T])$, which serves as the guessed propagation for an agent's decision-making, the optimal stopping problem for an agent starting at time $t$ and state $x$ is to choose a stopping time before the process first hits the boundary $\partial\Omega$. The corresponding value function is defined as
\[
u(x,t)=\min_\tau \mathbb{E}\left[\int_{t}^{\tau\wedge \tau_b} f(X_s,s,\mu)\,ds + \psi(X_{\tau\wedge\tau_b},\tau\wedge\tau_b,\mu)\,\middle|\, X_t=x\right],
\]
where $\tau_b = \inf\{s \ge t : X_s \in \partial\Omega\}$ is the first exit time of the process $X$ from the domain $\Omega$ after time $t$. We assume that $\psi=0$ when $x\in\partial \Omega$. Then
 by the dynamic programming process, $u$ satisfies the following PDE system:
\begin{equation}\label{equation_for_u}
\begin{cases}
\max(-\partial _tu-\mathcal{L} u-f(\mu),u-\psi(\mu))=0,&(x,t)\in \Omega\times(0,T);\\
u(x,T)=\psi(x,T,\mu),&x\in \Omega;\\
u(x,t)=0,&x\in \partial \Omega.
\end{cases}
\end{equation}
Similar as in the single agent problem, an optimal control of whether to leave can be implied from (\ref{equation_for_u}). On the other hand, given an optimal control taken by all agents, we can solve a Fokker-Planck equation that will be discussed below to obtain the actual crowd distribution $m$. An equilibrium is reached when $m=\mu$.

As discussed in the single-agent problem, given a guessed distribution $\mu \in \mathcal{M}(\Omega\times[0,T])$, the optimal control for an agent need not be unique: in some region of $\{u=\psi\}$, the agent may be indifferent between leaving and staying. If we require all agents to leave whenever $u=\psi$, we obtain the \textit{pure-strategy equilibrium}, which is restrictive and may fail to exist. If instead we allow each agent to choose a probability of leaving in the indifference region of $\{u=\psi\}$, we obtain the \textit{mixed-strategy equilibrium}, which is generally more satisfactory, as existence is guaranteed under mild regularity assumptions. The details of these two equilibrium notions are discussed in the next section.
\subsection{Equilibrium for OSMFG} 
\subsubsection{Pure Strategy Equilibrium}
For the pure strategy equilibrium, we require in addition that all agents choose to exit when $(x,t)\in\{u=\psi\}$ and stay when $(x,t)\in\{u<\psi\}$. The actual crowd propagation will satisfy the following Fokker-Planck equation: 
\begin{equation}\label{equation_for_m} 
\begin{cases}
\partial_tm-\mathcal{L}^* m=0,&(x,t)\in\{u<\psi\};\\
m=0,&(x,t)\in\{u=\psi\};\\
m(x,0)=m^0(x),&x\in\Omega;\\
m(x,t)=0,&x\in\partial\Omega.
\end{cases}
\end{equation}

When the actual crowd propagation matches the guessed one (i.e. $m=\mu$), we can combine (\ref{equation_for_u}) and (\ref{equation_for_m}) to obtain the PDE system of the Nash equilibrium with pure strategy:
\begin{equation}\label{pure_stategy}
\begin{cases}
\max(-\partial _tu-\mathcal{L} u-f(m),u-\psi(m))=0,&(x,t)\in \Omega\times(0,T);\\
\partial_tm-\mathcal{L}^* m=0,&(x,t)\in\{u<\psi\};\\
m=0,&(x,t)\in\{u=\psi\};\\
u(x,T)=\psi(x,T,m),\,m(x,0)=m^0(x),&x\in \Omega;\\
u(x,t)=0,\,m(x,t)=0,&x\in \partial \Omega.
\end{cases}
\end{equation}
We have mentioned that in some regions of the set $\{u=\psi\}$, agents are indifferent between exiting the game and staying in the game. And we have required all agents to leave in such cases. In other words, when an agent faces choices that bring the same cost, they should all choose a specific control rather than individually make a random choice. Therefore, we refer to (\ref{pure_stategy}) as the pure strategy equilibrium of mean field games with optimal stopping.

However, the pure strategy equilibrium may not exist for mean field game models. In fact, we have the following result:
\begin{theorem}\label{nonexist_pure}
	There exists $f$ and $\psi$ such that there is no solution for (\ref{pure_stategy}).
\end{theorem}

We refer to \cite{osmfgpde} for the proof of such a non-existence result for the stationary case. The construction of $f$ and $\psi$ for the time-dependent case  (\ref{pure_stategy}) is similar and we omit the details here. Instead of elaborating on the proof, we provide an economic intuition for such non-existence.

Consider a game with a population of agents facing an optimal stopping problem, where agents experience crowd aversion. As we will show, there exists no pure strategy Nash equilibrium in this setting. Suppose all other agents follow the pure strategy of exiting at time $t$. For any individual agent, the cost of staying past $t$ will be low due to crowd aversion. Thus, the agent's best response is to deviate from the crowd and continue playing. On the other hand, if all agents use the pure strategy of staying past $t$, an individual agent will prefer exiting at $t$ to avoid the crowd. In either case, no single pure strategy can be sustained in equilibrium, as agents have an incentive to deviate.

The underlying issue is that the definition of pure strategy equilibrium requires all agents to make exactly the same deterministic decision, even when indifferent. If we relax this and allow agents to make random decisions between leaving and staying in the game, an equilibrium may exist. This leads to the concept of the mixed strategy equilibrium, where agents follow randomized policies that make opponents indifferent across actions. By permitting randomization, mixed strategy equilibria are able to capture equilibrium behavior even when the pure strategy equilibria do not exist. This provides an avenue to model strategic interactions in mean field games with crowd aversion and optimal stopping.

\subsubsection{Mixed Strategy Equilibrium}

If we relax the constraint that forces agents to uniformly exit at $\{u=\psi\}$, the mixed strategy equilibrium can be defined. With multiple optimal choices available, mixed equilibria allow agents to randomize over strategies freely, rather than requiring deterministic decisions.

Assume that the cost functions $f$ and $\psi$ depend only on the crowd distribution $m$. Then the PDE system governing the mixed strategy equilibrium for OSMFG  is as follows:	\begin{align}\label{mixed_strategy}
\begin{cases}
\max(-\partial_tu-\mathcal{L} u-f(m),u-\psi(m))=0, &(x,t)\in\Omega\times(0,T);\\
\partial_tm-\mathcal{L}^* m=0, &(x,t)\in  \{u<\psi\};\\
\partial_tm-\mathcal{L}^* m\le0,\,m\ge 0, &(x,t)\in  \Omega\times(0,T);\\
\int_{\{u=\psi\}}(f(m)+(\partial_t+\mathcal{L})\psi(m))m\,dxdt=0;\\
u(x,T)=\psi(x,T,m),\,m(x,0)\le m^0(x),&x\in\Omega;\\
u=0,\,m=0,&(x,t)\in \partial \Omega\times (0,T);
\end{cases}
\end{align}
We refer to \cite{osmfgpde} for more discussions of the mixed strategy equilibrium. The complementary condition of (\ref{mixed_strategy}),
$$
\int_{\{u=\psi\}}(f(m)+(\partial_t+\mathcal{L})\psi(m))m\,dxdt=0,
$$
shows that $$\{u=\psi\}\cap \{f(m)+(\partial_t+\mathcal{L})\psi(m)=0\}$$ is the region where agents can choose to exit or to stay. It can be derived by standard stochastic calculus that agents would have the same cost whenever exiting or staying at $(x,t)$ when $(x,t)\in\{u=\psi\}\cap\{f(m)+(\partial_t+\mathcal{L})\psi(m)=0\}$.  
This is the main relaxation compared to the pure strategy equilibrium. 
\subsection{Existence of the Mixed Strategy Equilibrium}
The following existence result shows that the mixed-strategy equilibrium always exists under mild conditions, whereas the pure-strategy one may not. This makes the former a more suitable notion for OSMFGs.  In order to state the theorem rigorously, we first introduce several spaces and some technical assumptions for $\Omega$, $f$, $\psi$ and $m^0$.
\begin{definition}\label{spaceABC}
    We define function spaces $\mathcal{A}$, $\mathcal{B}$ and $\mathcal{C}$ as:
    $$\mathcal{A}=L^2\left((0, T), H_0^1(\Omega) \cap H^2(\Omega)\right) \cap H^1\left((0, T), L^2(\Omega)\right),$$
    $$\mathcal{B}=L^2\left((0, T), H_0^1(\Omega)\right),\quad\mathcal{C}=L^2\left((0, T), L^2(\Omega)\right).$$
\end{definition}

\begin{assumption}\label{assump_f&psi}
\begin{itemize}
\item $\Omega$ is a connected open bounded subset of $\mathbb{R}^d$ with a smooth boundary.
\item $m^0\in L^{\infty}(\Omega),$ and $m^0\ge 0$.
\item The running cost $f(\cdot,\cdot,m)$ and  the stopping cost $\psi(\cdot,\cdot,m)$ satisfy:
\begin{enumerate}
    \item The map $m\mapsto f(\cdot,\cdot,m)$ is continuous from $\mathcal{C}$ to itself;

    \item  The map $m\mapsto \psi(\cdot,\cdot,m)$ is continuous from $\mathcal{C}$ to $\mathcal{A}$.
\end{enumerate}
\end{itemize}
\end{assumption}

A pair $(u, m) \in \mathcal{A} \times \mathcal{B}$ is called a solution to the mixed strategy equilibrium system \eqref{mixed_strategy}, if $u$ and $m$ satisfy \eqref{mixed_strategy} in the sense of measure, as defined in \cite{osmfgpde}. More precisely, we introduce the following definition:
\begin{definition}\label{def:mix_solu}
    We say a pair $(u,m)\in \mathcal{A}\times \mathcal{B}$ is a solution to the mixed strategy equilibrium system \eqref{mixed_strategy}, if $(u,m)$ satisfies:
    \begin{itemize}
        \item $\max(-\partial_tu-\mathcal{L} u-f(m),u-\psi(m))=0,\quad\text{in }\mathcal{D}'(\Omega\times (0,T))$;
        \item $\partial_tm-\mathcal{L}^* m=0,\quad\text{in }\mathcal{D}'(\{u<\psi\})$;
        \item $\partial_tm-\mathcal{L}^* m\le0,\quad \text{in }\mathcal{D}'(\Omega\times(0,T))$, and $m\ge 0,\quad a.e.$;
        \item $\int_{\{u=\psi\}}(f(m)+(\partial_t+\mathcal{L})\psi(m))m\,dxdt=0$;
        \item $u(x,T)=\psi(x,T,m)$; and for any smooth $v$ such that $v\le 0$, $v(T)=0$,\footnote{The following condition \eqref{mix:initial_cond} actually implies that $m|_{t=0}\le m^0$.}
        \begin{equation}\label{mix:initial_cond}
        \int_{0}^{T}\int_{\Omega}\left(-\partial_tv-\Delta v\right)m-\int_{\Omega}m^0v(0)\ge 0.
        \end{equation}
    \end{itemize}
\end{definition}
Then the existence result can be formulated in the following theorem. 
\begin{theorem}\label{exist_mixed}
	Assume that $\mathcal{L}=\Delta.$ Suppose that Assumption \ref{assump_f&psi} holds. Then there exists at least one solution $(u,m)\in \mathcal{A}\times \mathcal{B}$ to (\ref{mixed_strategy}).
\end{theorem}

We note that a proof for this theorem was presented in Theorem 2.1 in \cite{osmfgpde}. However, as pointed out in \cite{measureflow}, the proof in\cite{osmfgpde} is flawed. To be precise, the weak convergence of \( m^\epsilon \) in \( \mathcal{B} \) does not guarantee the convergence of \( \int f(m^\epsilon) \, dm^\epsilon \) for nonlinear \( f(m^\epsilon) \), contrary to the claim in \cite{osmfgpde}. 

To establish the convergence of \( \int f(m^\epsilon) \, dm^\epsilon \), it is essential to demonstrate the strong convergence of \( m^\epsilon \) in the space \( \mathcal{C} \). Typically, strong convergence in $\mathcal{C}$ follows from weak convergence in $\mathcal{B}$ combined with suitable regularity of the time derivative, as provided by the Aubin--Lions lemma. In the present problem, however, such regularity is not available. Nevertheless, within the mixed strategy equilibrium setting \eqref{mixed_strategy}, we restrict to functions $m \in \mathcal{B}$ satisfying $\partial_t m - \Delta m \le 0$. This monotonicity property ensures strong sequential compactness: whenever $\{m^\epsilon\}$ is uniformly bounded in $\mathcal{B}$, a subsequence converges strongly in $\mathcal{C}$. The following lemma formalizes this result, which plays a pivotal role in the proof of the existence theorem.
\begin{lemma}\label{lem:weak_strong_comp}
    We define the set $\mathcal{T}$ as follows:
    \begin{equation}\label{def:mathcalT}
        \mathcal{T}=\{m\in \mathcal{B}|m\ge 0,\,(\partial_t-\Delta)m\le 0,\,m|_{t=0}\le 
        m^0\},
    \end{equation}
    where $m\in\mathcal{T}$ satisfies all the conditions (i.e. $m\ge 0,\,(\partial_t-\Delta)m\le 0,\,m|_{t=0}\le m^0$) in the sense of Definition \ref{def:mix_solu}. Assume that $\{m_n\}_{n=1}^{\infty}\subset\mathcal{T}$ is bounded in $\mathcal{B}$, then there exists a subsequence $\{m_{k_n}\}_{n=1}^{\infty}\subset\{m_{n}\}_{n=1}^{\infty}$, and $m^*\in \mathcal{C}$, s.t. $m_{k_n}\rightarrow m^*(n\rightarrow\infty)$ in the sense of $\|\cdot\|_{\mathcal{C}}$.
\end{lemma}
\begin{proof}
    Suppose that $\|m_n\|_{\mathcal{B}}\le M$ for all $n$. The proof consists of two steps. In the first step, we show that $\{m_n\}$ is $L^1$-continuous in both the time-variable and the spatial-variable.  An application of the Kolmogorov-Riesz-Fr\'echet theorem then yields a strongly convergent subsequence in $L^1(\Omega\times (0,T))$. In the second step, we show that we can extract a further subsequence that  converges strongly in $\mathcal{C}$ from the $L^1$-convergent subsequence obtained from the first step. 
    \begin{itemize}
        \item \textbf{Step 1.} We first show that $\{m_n\}_{n=1}^{\infty}$  is $L^1$-equicontinuous in the time-variable. Indeed, for any $m\in \mathcal{T}$ such that $\|m\|_{\mathcal{B}}\le M$, and any Lebesgue points $t_1<t_2$\footnote{It follows from the Bochner version of the Lebesgue differentiation theorem that for $m \in L^p((0,T); X)$, where $X$ is a separable Banach space, almost every $t \in (0,T)$ is a strong Lebesgue point, i.e., $\lim_{h\to 0^+} \frac{1}{h} \int_{t}^{t+h} \|m(s) - m(t)\|_X \, ds = 0$.}, we have
        $$
        \begin{aligned}
        \|m(\cdot,t_1)-m(\cdot,t_2)\|_{L^1(\Omega)}&=\int_{\Omega}|m(x,t_1)-m(x,t_2)|dx\\
        &=\int_{\Omega}m(x,t_1)-m(x,t_2)dx+2\int_{\Omega}\left[m(x,t_2)-m(x,t_1)\right]_+dx\\
        &\le M(t_1)-M(t_2)+2\int_{\Omega}\left[e^{(t_2-t_1)\Delta}m(t_1)-m(t_1)\right]_+dx,
        \end{aligned}
        $$
        where $M(t)$ denotes the total mass, which is decreasing in $t$ since $\partial_tm-\Delta m\le 0$. In the last line, we used the fact that $m(\cdot,t_2)\le e^{(t_2-t_1)\Delta}m(\cdot,t_1)$. Hence, for $\tau>0$ and $t\in [0,T-\tau]$,
        $$
        \begin{aligned}
        \|m(\cdot,t)-m(\cdot,t+\tau)\|_{L^1(\Omega)}&\le M(t)-M(t+\tau)+2\int_{\Omega}\left[e^{\tau\Delta}m(x,t)-m(x,t)\right]_+dx\\
        &\le M(t)-M(t+\tau)+2|\Omega|^{\frac{1}{2}}\|e^{\tau\Delta}m(\cdot,t)-m(\cdot,t)\|_{L^2(\Omega)}\\
        &\le M(t)-M(t+\tau)+C|\Omega|^{\frac{1}{2}}\tau^{\frac{1}{2}}\|\nabla m(\cdot,t)\|_{L^2(\Omega)}.
        \end{aligned}
        $$
        In the last line, we used the spectral property of the heat kernel on $\Omega$ with homogeneous Dirichlet boundary condition on $\partial\Omega$. Therefore, 
        $$
        \begin{aligned}
        &\int_{0}^{T-\tau}\|m(\cdot,t)-m(\cdot,t+\tau)\|_{L^1(\Omega)}dx\\
        \le &\int_{0}^{T-\tau}M(t)dt-\int_{\tau}^{T}M(t)dt+C|\Omega|^{\frac{1}{2}}\tau^{\frac{1}{2}}\int_{0}^{T-\tau}\|\nabla m(\cdot,t)\|_{L^2(\Omega)}dt\\
        \le &\int_{0}^{\tau}M(t)dt-\int_{T-\tau}^{T}M(t)dt+C|\Omega|^{\frac{1}{2}}\tau^{\frac{1}{2}}T^{\frac{1}{2}}\left(\int_{0}^{T}\|\nabla m(\cdot,t)\|^2_{L^2(\Omega)}dt\right)^{\frac{1}{2}}\\
        \le &\tau M(0)+C|\Omega|^{\frac{1}{2}}\tau^{\frac{1}{2}}T^{\frac{1}{2}}\left(\int_{0}^{T}\|\nabla m(\cdot,t)\|^2_{L^2(\Omega)}dt\right)^{\frac{1}{2}}.
        \end{aligned}
        $$
        This implies that $m$ is $L^1$-equicontinuous in the time-variable. The equicontinuity in the spatial-variable is a direct consequence of the uniform bound $\|m\|_{L_T^2\dot{H}^{1}(\Omega)}\le \|m\|_{\mathcal{B}}\le M$. According to Kolmogorov-Riesz-Fr\'echet Theorem(Theorem 4.26 in \cite{brezis2011functional}), there exists a subsequence of $\{m_n\}_{n=1}^{\infty}$, denoted by $\{m_{\ell_n}\}_{n=1}^{\infty}$, such that $m_{\ell_n}$ converges to some $m^*$ strongly in $L^1(\Omega\times(0,T))$.
        \item \textbf{Step 2. }According to Riesz's Theorem, there exists a subsequence $\{m_{k_n}\}_{n=1}^{\infty}\subset \{m_{\ell_n}\}_{n=1}^{\infty}$, such that $m_{k_n}\rightarrow m^*,\,a.e.$. Since $\Omega$ is bounded, and $m_{k_n}\le e^{t\Delta} m^0\, a.e.$, we conclude that $m_{k_n}\rightarrow m^*$ strongly in $L^2(\Omega\times (0,T))$, and hence, in $\mathcal{C}$. Therefore we finish the proof.
    \end{itemize}
\end{proof}

With Lemma \ref{lem:weak_strong_comp}, the proof for Theorem \ref{exist_mixed} can proceeded as in the proof for Theorem 2.1 in \cite{osmfgpde}, and we omit the details here. 

Although the existence of the mixed strategy equilibrium can be guaranteed, it is yet difficult to design a simple and efficient numerical PDE algorithm for (\ref{mixed_strategy}) directly due to its complicated form. The coupling between $u$ and $m$ is more intricate  due to the complementary condition  
$$\int_{\{u=\psi\}}(f(m)+(\partial_t+\mathcal{L})\psi(m))mdxdt=0 .$$ 
In the next section, the primary goal is to construct an iterative algorithm that provides an acccurate approximation of the mixed strategy equilibria for OSMFGs via a generalized fictitious play.

\section{Algorithm Construction}
In this section, we introduce an iterative algorithm for solving the mixed-strategy equilibrium system \eqref{mixed_strategy}. The key idea is to obtain the mixed-strategy equilibrium by applying a generalized fictitious play to the more tractable pure-strategy system \eqref{pure_stategy}. We proceed as follows. We first introduce the generalized fictitious play and prove its convergence from pure-strategy to mixed-strategy equilibria for potential games. We then present finite difference schemes for the obstacle and Fokker--Planck equations. Finally, we show that the fully implicit scheme preserves this convergence property in the discrete setting.

\subsection{Fictitious Play}

Fictitious play is a learning procedure in which, at each iteration step, every agent chooses the best response strategy based on the average of the strategies in the previous iterations employed by other agents. It is therefore natural to expect that the best response strategy converges to a Nash equilibrium of the game. For classical MFG models, it has been shown in \cite{fictitiousplay} that the sequence obtained via fictitious play converges to the MFG solution for potential games. For OSMFG models, \cite{LP} proves that the sequence obtained via the fictitious-play-based linear programming method converges to the relaxed Nash equilibrium in the measure flow framework for potential games. However, fictitious play has not yet been directly applied to the PDE systems of Nash equilibrium for OSMFG models. 

Our aim is to adapt fictitious play to the PDE systems in order to compute mixed-strategy equilibrium solutions for OSMFG models.  This requires novel theoretical developments, as the existing convergence results for fictitious play in classical MFG settings \cite{fictitiousplay} do not involve mixed-strategy equilibria and hence do not directly extend to OSMFGs.

Recall that fictitious play approximates a mixed-strategy equilibrium through iterative pure strategies. In a pure strategy, an agent decides deterministically whether to quit; in a mixed strategy, he specifies a quitting probability. This is the key distinction between the two equilibrium notions. In fictitious play for OSMFG, at each iteration, every agent chooses a pure strategy that is the best response to the average historical distribution. Over repeated plays, the agent may adopt different pure strategies at the same state $(x,t)$ in different iterations---quitting in some iterations and staying in others. Consequently, his overall quitting behavior at $(x,t)$ is described by a probability between 0 and 1, which illustrates the idea behind a mixed strategy.

Each iteration of a fictitious play for OSMFG can be divided into the following three steps:
\begin{enumerate}
\item{\bf{Finding the Best Response.}} Agents calculate an optimal pure strategy response based on the so-called updated distribution of agents, which encodes the historical information from previous iterations.

\item{\bf{Propagating the Proposed Distribution.}} The proposed distribution of agents is determined based on the strategies obtained in the previous step.

\item{\bf{Calculating the Updated  Distribution.}} The agents calculate the updated distribution by incorporating the proposed distribution with a certain updating rule. The updated distribution will be used to calculate controls in the next iteration.
\end{enumerate}
The procedure can be formulated in the definition below.

\begin{definition}\label{fictitiousplay}
	(generalized fictitious play) Given the initial distribution $m^0$, and the initial averaged distribution $\bar{m}_0$, the following iteration is called a generalized fictitious play.
	\begin{equation}\label{fictitious}
		\begin{cases}
		\max(-\partial _tu_{n+1}-\mathcal{L} u_{n+1}-f(\bar{m}_n),u_{n+1}-\psi(\bar{m}_n))=0,&(x,t)\in \Omega\times(0,T);\\
		\partial_tm_{n+1}-\mathcal{L}^* m_{n+1}=0,&(x,t)\in\{u_{n+1}<\psi(\bar{m}_n)\};\\
		m_{n+1}=0,&(x,t)\in\{u_{n+1}=\psi(\bar{m}_n)\};\\
		u_{n+1}(x,T)=\psi(x,T,\bar{m}_n),\,m_{n+1}(x,0)\le m^0(x),&x\in \Omega;\\
		u_{n+1}=0,\,m_{n+1}=0,&x\in \partial \Omega;\\
		\end{cases}
	\end{equation}
	where $n=0,1,2,...$ is the iteration round, and
 \begin{equation}
     \bar{m}_n=\delta_n m_n+(1-\delta_n)\bar{m}_{n-1}.
 \end{equation}
 Here, $m_n$ is called the proposed distribution, and $\bar{m}_n$ is called the updated distribution (or averaged distribution). $\delta_n$ is the learning rate satisfying
 \begin{equation} \label{cond_delta}
\sum\delta_n\rightarrow\infty,\quad\sum\delta_n^2<\infty.
 \end{equation}
\end{definition}

With $\delta_n=1/n$, the generalized fictitious play reduces to the classic fictitious play. Generalized Fictitious play can be viewed as a learning process for the agents. In each iteration $n$, agents have observed all past proposed distributions $m_0,...,m_n$ and calculate the updated distribution $\bar{m}_n$. The updated distribution $\bar{m}_n$ is a weighted average and is taken as the guessed distribution for the next iteration. Using $\bar{m}_n$, agents then obtain a new pure strategy response, yielding another proposed distribution $m_{n+1}$. This distribution $m_{n+1}$ provides new information such that agents can obtain the updated distribution $\bar{m}_{n+1}$ for the next iteration by averaging $m_0$, $\cdots$, $m_n$, $m_{n+1}$. In other words, agents consecutively update the updated distribution $\bar{m}_{n+1}$ by including the latest proposed distribution $m_{n+1}$. 

Note that for each iteration, $m_n$ and $\bar{m}_n$ represent both the crowd distribution and an individual agent's state distribution. Specifically, they characterize a pure strategy and a mixed strategy for a representative agent. As $n\rightarrow\infty$, we expect the sequence of updated distribution $\bar{m}_n$ to converge to the distribution in the mixed strategy equilibrium. We will prove this convergence when the game is a potential game in the next section.

\subsection{Convergence of Fictitious Play to Mixed Strategy Equilibrium}

For simplicity, we assume the generator of the diffusion process $\mathcal{L}$ is simply the Laplacian operator $\Delta$ in the subsequent analysis. To prove the convergence of the generalized fictitious play (\ref{fictitious}), we first need to ensure that the iteration scheme is well-defined, i.e. $(u_n,\bar{m}_n)$ should belongs to $\mathcal{A}\times \mathcal{B}$ for all $n\in\mathbb{N}_+$. We have the following theorem.
\begin{theorem}\label{thm:regular_point}
   Assume that $\mathcal{L}=\Delta$. We suppose that $\bar{m}_n\in \mathcal{B}$ and Assumption \ref{assump_f&psi} holds. Then there exists $(u_{n+1},m_{n+1})\in \mathcal{A}\times \mathcal{B}$ such that $(u_{n+1},m_{n+1})$ is the solution to system  \eqref{fictitious}. 
\end{theorem}
\begin{proof}
     We notice that equation \eqref{fictitious} can be separated into an backward equation for $u_{n+1}$ and forward equation for $m_{n+1}$. The proof for this theorem can then be broken into two parts. To be precise, first, given $\bar{m}_n\in \mathcal{B}$, we prove that there exists $u_{n+1}\in\mathcal{A}$, such that $u_{n+1}$ is the solution to the following system:
     \begin{equation}\label{eq:u_fictitious}
       \begin{cases}
		\max(-\partial _tu_{n+1}-\Delta u_{n+1}-f(\bar{m}_n),u_{n+1}-\psi(\bar{m}_n))=0,&(x,t)\in \Omega\times(0,T);\\
  u_{n+1}(x,T)=\psi(x,T,\bar{m}_n),&x\in\Omega;\\
  u_{n+1}=0,&x\in \partial \Omega.
        \end{cases}
     \end{equation}
      Second, given $u_{n+1}\in\mathcal{A}$, we need to prove that there exists $m_{n+1}\in\mathcal{B}$, such that $m_{n+1}$ is the solution to the following system:
     \begin{equation}\label{eq:m_fictitious}
         \begin{cases}
             \partial_tm_{n+1}-\Delta m_{n+1}=0,&(x,t)\in\{u_{n+1}<\psi\};\\
  m_{n+1}=0,&(x,t)\in\{u_{n+1}=\psi\};\\
	m_{n+1}(x,0)\le m^0(x),&x\in \Omega;\\
		m_{n+1}=0,&x\in \partial \Omega.\\
         \end{cases}
     \end{equation}
     
     Since $f(\bar{m}_n)\in \mathcal{C}$ and $\psi(\bar{m}_n)\in \mathcal{A}$, we can directly conclude a solution $u_{n+1}\in \mathcal{A}$ to \eqref{eq:u_fictitious} according to the regularity results of obstacle equation, see Section 5 in \cite{osmfgpde} for details. 
     
     To finish this proof, it suffices to prove that there exists $m_{n+1}\in\mathcal{B}$, such that $m_{n+1}$ is the solution to \eqref{eq:m_fictitious}. We consider the following $\epsilon-$penalized system to system \eqref{eq:m_fictitious} as follows:
   \begin{equation}\label{eq:penalized_m}
   \begin{cases}
       	\partial_t m_\epsilon-\Delta m_{\epsilon}+\frac{1}{\epsilon}1_{\{u_{n+1}=\psi\}}m_{\epsilon}=0,& \text{in } \Omega\times(0,T];\\
		m_{\epsilon}(0,\cdot)=m^0(\cdot), &\text{in } \Omega;\\
        m_{\epsilon}=0, &x\in \partial \Omega.
  \end{cases}
   \end{equation}
   By the classic theory for linear parabolic equations, we know that there exists 
   $$m_{\epsilon}\in L^2(0,T;H_0^1(\Omega))\cap H^1(0,T;H^{-1}(\Omega)),$$ such that $m_{\epsilon}$ is a solution to system \eqref{eq:penalized_m}. The remainder of the proof can be divided into two main steps. In the first step, we show that we have an uniformly estimation of $\|m_{\epsilon}\|_{L^2(0,T;H_0^1(\Omega))}$, and hence there exists $m\in L^2(0,T;H_0^1(\Omega))$ such that $m_{\epsilon}\rightharpoonup m$. In the second step, we verify that $m$ is the solution to \eqref{eq:m_fictitious}.
\begin{enumerate}
\item We multiple $m_{\epsilon}$ to \eqref{eq:penalized_m} and obtain
\begin{equation}\label{multi}
m_{\epsilon}\partial_t m_{\epsilon}-m_{\epsilon}\Delta m_\epsilon +\frac{1}{\epsilon}1_{\{u=\psi\}}m_{\epsilon}^2=0.
\end{equation}
We integral (\ref{multi}) in $[0,t]\times  \Omega$ for any $t\in[0,T]$, and have the first estimate
\begin{equation}\label{estimate1}
	\frac{1}{2}\left(||m_\epsilon(t,\cdot)||_{L^2(\Omega)}^2-||m_\epsilon(0,\cdot)||_{L^2(\Omega)}^2\right)=-\int_{0}^{t}||\nabla m_{\epsilon}(\tau,\cdot)||_{L^2(\Omega)}^2d\tau-\frac{1}{\epsilon}\int_{\{u=\psi\}\cap\left(\Omega\times[0,t]\right)}m_{\epsilon}^2dxd\tau\le 0.
\end{equation}
We integral (\ref{multi}) in $[0,T]\times \partial \Omega$, and have the second estimate
\begin{equation}\label{estimate2}
	\frac{1}{2}\left(||m_\epsilon(T,\cdot)||_{L^2(\Omega)}^2-||m_\epsilon(0,\cdot)||_{L^2(\Omega)}^2\right)+\int_{0}^{T}||\nabla m_{\epsilon}(t,\cdot)||_{L^2(\Omega)}^2dt=-\frac{1}{\epsilon}\int_{\{u=\psi\}}m_{\epsilon}^2dxdt\le 0.
\end{equation}
From (\ref{estimate1}), we know that for any $t\in[0,T]$, we have
$$
||m_\epsilon(t,\cdot)||_{L^2(\Omega)}^2\le ||m_\epsilon(0,\cdot)||_{L^2(\Omega)}^2=||m^0||_{L^2(\Omega)}^2,
$$
yielding that
\begin{equation}\label{estimate_result1}
\int_{0}^{T}||m_{\epsilon}(t,\cdot)||_{L^2(\Omega)}^2dt\le T||m^0||_{L^2(\Omega)}^2.
\end{equation}
From (\ref{estimate2}), we know
\begin{equation}\label{estimate_result2}
	\int_{0}^{T}||\nabla m_{\epsilon}(t,\cdot)||_{L^2(\Omega)}^2dt\le \frac{1}{2}||m^0||_{L^2(\Omega)}^2.
\end{equation}
We combine (\ref{estimate_result1}) and (\ref{estimate_result2}) and obtain an uniform bound for $||m_{\epsilon}||_{L^2(0,T;H_{0}^1(\Omega))}^2$
\begin{equation}\label{unif_bound}
	||m_{\epsilon}||_{L^2(0,T;H_{0}^1(\Omega))}^2\le (T+\frac{1}{2})||m^0||_{L^2(\Omega)}^2.
\end{equation}
Thus we can construct a sequence $\{\epsilon_{n}\}\rightarrow 0$ and find $m\in L^2(0,T;H_{0}^1(\Omega))$, such that $
m_{\epsilon_n}\rightharpoonup m.
$
\item In order to verify $m$ is the solution to \eqref{eq:m_fictitious}, it suffices to verify
\begin{enumerate}
    \item\label{emua} $\partial_tm-\Delta m=0,\quad \text{in }\mathcal{D}'(\{u_{n+1}<\psi\})$; 
    \item\label{emub} $m=0$,\quad \text{a.e. in }$\{u_{n+1}=\psi\}$.
    \item \label{emuc} $m|_{t=0}\le m^0$, i.e. for any smooth $v$ such that $v\le 0$, $v(T)=0$,
        \begin{equation}\label{pure:initial_cond}
        \int_{0}^{T}\int_{\Omega}\left(-\partial_tv-\Delta v\right)m-\int_{\Omega}m^0v(0)\ge 0.
        \end{equation}
\end{enumerate}

For (\ref{emua}), we consider $\phi\in C_{c}^{\infty}(\{u_{n+1}<\psi\})$. We multiply $\phi$ with \eqref{eq:penalized_m} and integral over $\{u_{n+1}<\psi\}$ to obtain
$$
\int_{\{u_{n+1}<\psi\}}(\phi\partial_tm_\epsilon-\phi\Delta m_{\epsilon})dxdt=\int_{\{u_{n+1}<\psi\}}(-m_{\epsilon}\partial_t\phi - m_{\epsilon}\Delta\phi) dxdt=0.
$$
By the weak convergence $m_{\epsilon_n}\rightharpoonup m$ in $L^2(0,T;H_{0}^1(\Omega))$, we know that
$$
\int_{\{u_{n+1}<\psi\}}(-m\partial_t\phi - m\Delta\phi) dxdt=0,
$$
and therefore
$$
\partial_tm-\Delta m=0,\quad\text{ in }\mathcal{D}'(\{u_{n+1}<\psi\}).
$$

For (\ref{emub}), according to \eqref{estimate2}, we have
$$
\begin{aligned}
\int_{\{u=\psi\}} m_{\epsilon}^2dxdt&=-\frac{1}{2}\epsilon \left(||m_\epsilon(T,\cdot)||_{L^2(\Omega)}^2-||m_\epsilon(0,\cdot)||_{L^2(\Omega)}^2\right)-\int_{0}^{T}||\nabla m_{\epsilon}(t,\cdot)||_{L^2(\Omega)}^2dt\\
&\le \frac{1}{2}\epsilon||m_\epsilon(0,\cdot)||_{L^2(\Omega)}^2\rightarrow 0.
\end{aligned}
$$
Hence by the Cauchy-Schwarz Inequality, we have
$$
\int_{\{u=\psi\}}m_{\epsilon}dxdt\le \sqrt{|\{u=\psi\}|}\sqrt{\int_{\{u=\psi\}}m_{\epsilon}^2dxdt}\rightarrow 0,\,\text{ as }\epsilon\rightarrow 0.
$$
We take the test function $\phi=1_{\{u=\psi\}}$. By the weak convergence of $m_{\epsilon_n}$, we have
$$
\int_{\{u=\psi\}}mdxdt=0.
$$
Hence $m=0$ a.e. in $\{u=\psi\}$.

For (\ref{emuc}), the inequality \eqref{pure:initial_cond} can be obtained by passing $\epsilon\rightarrow 0$ in
$$
\int_{0}^{T}\int_{\Omega}\left(-\partial_tv-\Delta v\right)m_{\epsilon}-\int_{\Omega}m^0v(0)\ge 0.
$$
Therefore, we finish the proof.
\end{enumerate}
\end{proof}
\begin{remark}
    It is straightforward to see that $m_{n+1}\in \mathcal{T}$ defined in \eqref{def:mathcalT}.
\end{remark}

In what follows, we impose some additional technical assumptions on the functions $f$ and $\psi$. A common assumption is that $f-(\partial_t-\Delta)\psi$ is the variation of a potential function. Games that satisfy this property are known as potential games, of which we present the formal definition as follows.

\begin{definition}{(potential games)}\label{potentialgames}
	We call an optimal stopping mean field game a potential game if there exist a potential function $F\in C^{1,1}( \mathcal{C})$ such that
	$$
	f+(\partial_t+\Delta)\psi=\frac{\delta F}{\delta m}.
	$$
\end{definition}
Now we present the main convergence result.

\begin{theorem}\label{fictitious_convergence}
 Suppose that Assumption \ref{assump_f&psi} holds and the game is a potential game as in Definition \ref{potentialgames}. Assume that $\mathcal{L}=\Delta$. For any $n\in \mathbb{N}$, let $(u_n,\bar{m}_n)$ be the $n^{th}$ iteration result generated by the fictitious play (\ref{fictitiousplay}). Then any cluster point (in the sense of $\|\cdot\|_{\mathcal{A}\times\mathcal{C}}$) of the sequences $(u_n,\bar{m}_n)$  is a solution to (\ref{mixed_strategy}). 
\end{theorem}  

\begin{proof}
	Since $\psi\in \mathcal{A}$, we have $(\partial_t-\Delta)\psi\in \mathcal{C}$. Without loss of generality, we assume $\psi=0$ (otherwise we consider $\tilde{u}:=u-\psi,\,\tilde{m}=m,\,\tilde{f}=f+(\partial_t+\Delta)\psi$ and $\tilde{\psi}=0)$.

	The rest of the proof is divided into 3 steps. In the first step, we aim to prove that for any $n\in \mathbb{N}$ and $m\in\mathcal{T}$, 
 \begin{equation}\label{step1_result}
 \left<\frac{\delta F(\bar{m}_n)}{\delta m},m_{n+1}\right>\le \left<\frac{\delta F(\bar{m}_n)}{\delta m},m\right>,
 \end{equation}
 where $<\cdot,\cdot>$ denotes the inner product in $\mathcal{C}$. In the second step, we utilize (\ref{step1_result}) to prove that any regular cluster point $(u^*,m^*)$ should satisfy that $m^*$ is a local minimizer of $F$. In the third step, we show that any regular cluster point $(u^*,m^*)$ is a solution to (\ref{mixed_strategy}).
	\begin{enumerate}
	\item 	Recall that we define the test set $\mathcal{T}$ as follows:
    $$
    \mathcal{T}=\{m\in \mathcal{B}|m\ge 0, (\partial_t-\Delta)m\le 0, m|_{t=0}\le m^0\}.
    $$
In this step, we want to prove that, for any $n\in \mathbb{N}$ and $m\in \mathcal{T}$,
	\begin{equation}\label{mid1theo3.1}
	\left<\frac{\delta F(\bar{m}_n)}{\delta m},m_{n+1}-m\right>=\int_{0}^{T}\int_{\Omega}f(\bar{m}_n)(m_{n+1}-m)dxdt\le 0.
	\end{equation}
 
Denote $\Omega_{n+1}^1=\{u_{n+1}<0\},\,\Omega_{n+1}^2=\{u_{n+1}=0\}$. Then $\Omega\times(0,T)=\Omega_{n+1}^1\cup\Omega_{n+1}^2$.

In $\Omega_{n+1}^1$, we have $-\partial_tu_{n+1}-\Delta u_{n+1}=f(\bar{m}_n)$ from the first equation of (\ref{fictitiousplay}). Therefore, in $\Omega_{n+1}^1$, for any $m\in\mathcal{T}$, we have
\begin{equation}\label{step1_Omega1}
(-\partial_tu_{n+1}-\Delta u_{n+1}-f(\bar{m}_n))(m_{n+1}-m)=0.
\end{equation}

In $\Omega_{n+1}^2$, we have 
\begin{equation}\label{step1_Omega2_u}
-\partial_tu_{n+1}-\Delta u_{n+1}\le f(\bar{m}_n)
\end{equation}
from the first equation of (\ref{fictitiousplay}). Additionally, we have $m_{n+1}=0$ in $\Omega_{n+1}^2$ from the third equation of (\ref{fictitiousplay}). For any $m\in\mathcal{T}$, $m\ge 0$ and hence in $\Omega_{n+1}^2$ 
\begin{equation}\label{step1_Omega2_m}
m_{n+1}-m\le 0.
\end{equation}
Combining (\ref{step1_Omega2_u}) and (\ref{step1_Omega2_m}), in $\Omega_{n+1}^2$, for any $m\in\mathcal{T}$, we have
\begin{equation}\label{step1_Omega2}
(-\partial_tu_{n+1}-\Delta u_{n+1}-f(\bar{m}_n))(m_{n+1}-m)\ge0.
\end{equation}

(\ref{step1_Omega1}) and (\ref{step1_Omega2}) yield that 
	\begin{equation}\label{mid2theo3.1}
	\begin{aligned}
	&\int_{0}^{T}\int_{\Omega}f(\bar{m}_n)(m_{n+1}-m)dxdt\\
    =&\int_{\Omega_{n+1}^1}f(\bar{m}_n)(m_{n+1}-m)dxdt+\int_{\Omega_{n+1}^2}f(\bar{m}_n)(m_{n+1}-m)dxdt\\
	\le &\int_{\Omega_{n+1}^1}(-\partial_tu_{n+1}-\Delta u_{n+1})(m_{n+1}-m)dxdt+\int_{\Omega_{n+1}^2}(-\partial_tu_{n+1}-\Delta u_{n+1})(m_{n+1}-m)dxdt\\
    = &\int_{\Omega_{n+1}^1}(-\partial_tu_{n+1}-\Delta u_{n+1})(m_{n+1}-m)dxdt,
	\end{aligned}
	\end{equation}
 where the last equality holds because
 $$
 -\partial_tu_{n+1}-\Delta u_{n+1}=0,\quad a.e. \text{ in }\Omega_{n+1}^2
 $$
 for $u_{n+1}\in \mathcal{A}$. The inequality \eqref{mid2theo3.1} holds as an equality if and only if (recall  $m_{n+1}=0$ in $\Omega_{n+1}^2$)
 \begin{equation}\label{cond1}
 \int_{\Omega_{n+1}^2}f(\bar{m}_n)m\,dxdt=0.
 \end{equation}
	Thus, to prove  (\ref{mid1theo3.1}), it suffices to show that the right hand side of (\ref{mid2theo3.1}) is not greater than 0. Indeed, it is straightforward to verify from \eqref{eq:penalized_m} that 
    $$
    \begin{aligned}
    &\int_{\Omega_{n+1}^1}(-\partial_t u_{n+1}-\Delta u_{n+1})m^\epsilon dxdt-\int_{\Omega}m^0(x)u(0,x)dx\\=&\int_{0}^T\int_{\Omega}(-\partial_tu_{n+1}-\Delta u_{n+1})m^{\epsilon}dxdt-\int_{\Omega}m^0(x)u(0,x)dx\\
    =&\int_{0}^T\int_{\Omega} u_{n+1} \cdot \frac{1}{\epsilon}1_{\{u_{n+1}=0\}}dx=0,
    \end{aligned}
    $$
    which implies that
    \begin{equation}\label{eq:mid1}
    \int_{\Omega_{n+1}^1}(-\partial_t u_{n+1}-\Delta u_{n+1})m_{n+1} dxdt-\int_{\Omega}m^0(x)u(0,x)dx=0,
    \end{equation}
    when taking $\epsilon\rightarrow0$.
    Since $m\in \mathcal{T}$, the condition $m|_{t=0}\le m^0$ yields that
 \begin{equation}\label{eq:mid2}
 \begin{aligned}
 &\int_{\Omega_{n+1}^1}(-\partial_t u_{n+1}-\Delta u_{n+1})mdxdt-\int_{\Omega}m^0u(0,x)dx\\
 =&\int_{0}^T\int_{\Omega}(-\partial_t u_{n+1}-\Delta u_{n+1})mdxdt-\int_{\Omega}m^0u(0,x)dx\ge 0.
 \end{aligned}
 \end{equation}
It then follows from \eqref{eq:mid1} and \eqref{eq:mid2} that
 $$
 \int_{\Omega_{n+1}^1}(-\partial_tu_{n+1}-\Delta u_{n+1})(m_{n+1}-m)dxdt\le 0.
 $$
The inequality holds as an equality if and only if \eqref{eq:mid2} holds as an equality, which implies that
	\begin{equation}\label{cond2}
	     \partial_t m-\Delta m=0\text{ in }\Omega_{n+1}^1\text{ in the distributional sense.}
	\end{equation}
 
	\item Let $(u^*,m^*)$ be a cluster point of $(u_n,\bar{m}_n)$, we claim that $m^*$ satisfies $$\left<\frac{\delta F}{\delta m}(m^*),m-m^*\right>\ge 0, \quad  \forall \, m\in\mathcal{T}.$$ If not, we can find $\bar{m}\in\mathcal{T}$ and $\left<\frac{\delta F}{\delta m}(m^*),\bar{m}-m^*\right>< 0$. Suppose that $(u^{**},m^{**})$ is the next iteration result of $(u^*,m^*)$ in the fictitious play (\ref{fictitiousplay}), i.e. the solution of the following system:
\begin{equation}\label{mstarstar}
	\begin{cases}
	\max(-\partial _tu^{**}-\Delta u^{**}-f(m^*),u^{**})=0,&(x,t)\in \Omega\times(0,T);\\
	\partial_tm^{**}-\Delta m^{**}=0,&(x,t)\in\{u^{**}<0\};\\
	m^{**}=0,&(x,t)\in\{u^{**}=0\};\\
	u^{**}(x,T)=0,\,m^{**}(x,0)=m^0(x),&x\in \Omega;\\
	u^{**}=0,\,m^{**}=0,&x\in \partial \Omega.
	\end{cases}
	\end{equation}
	From step 1 we know that
 $$
 \left<\frac{\delta F(m^*)}{\delta m},m^{**}\right>\le \left<\frac{\delta F(m^*)}{\delta m},\bar{m}\right>.
 $$
 The assumption 
 $\left<\frac{\delta F}{\delta m}(m^*),\bar{m}-m^*\right>< 0$ can be written as
	$$
	\left<\frac{\delta F(m^*)}{\delta m},\bar{m}\right>< \left<\frac{\delta F(m^*)}{\delta m},m^*\right>.
	$$ 
 Therefore, we have
 $$
	\left<\frac{\delta F(m^*)}{\delta m},m^{**}\right>< \left<\frac{\delta F(m^*)}{\delta m},m^*\right>.
	$$
 We define that
	$$
	\left<\frac{\delta F(m^*)}{\delta m},m^{**}-m^{*}\right>=\int_{0}^{T}\int_{\Omega}f(m^*)(m^{**}-m^*)dxdt:=-\ell<0.
	$$

Suppose that $\tilde{m}$ is the solution of the Fokker-Planck equation when no agents exit the game except through the boundary $\partial \Omega$:
	$$
	\begin{cases}
	\partial _tm-\Delta m=0,&(x,t)\in \Omega \times [0,T];\\
	m(x,0)=m^0(x),&x\in\Omega;\\
	m=0,&x\in\partial \Omega.
	\end{cases}
	$$
	 For any $m\in \mathcal{T}$, we have $\left|\left|m\right|\right|_{\mathcal{C}}\le \left|\left|\tilde{m}\right|\right|_{\mathcal{C}}$.
 Thus by the Lipschitz continuity of $f$ in $\mathcal{C}$, we conclude that $\left|\left|f(m)\right|\right|_{\mathcal{C}}$ is bounded for $m\in\mathcal{T}$, and there exists $\epsilon>0$, such that as long as $\left|\left|\bar{m}_n-m^*\right|\right|_{\mathcal{C}}<\epsilon$, we have
\begin{equation}\label{step2_conti1}
 \int_{0}^{T}\int_{\Omega}(f(\bar{m}_n)-f(m^*))(m^{**}-m^{*})dxdt<\frac{\ell}{4},
 \end{equation}
 and
 \begin{equation}\label{step2_conti2}
 \int_{0}^{T}\int_{\Omega}f(\bar{m}_n)(m^{*}-\bar{m}_n)dxdt<\frac{\ell}{4}.
 \end{equation}
 We can deduce from (\ref{step2_conti1}) and (\ref{step2_conti2}) that
 \begin{equation}\label{step2_continuous}
 \begin{aligned}
 &\int_{0}^{T}\int_{\Omega}f(\bar{m}_n)(m^{**}-\bar{m}_n)dxdt\\
 =&\int_{0}^{T}\int_{\Omega}f(\bar{m}_n)(m^{**}-m^{*})dxdt+\int_{0}^{T}\int_{\Omega}f(\bar{m}_n)(m^{*}-\bar{m}_n)dxdt\\
 <&\int_{0}^{T}\int_{\Omega}f(m^*)(m^{**}-m^{*})dxdt+\frac{\ell}{2}=-\frac{\ell}{2}.
 \end{aligned}
 \end{equation}

    It is obvious that $m^{**}\in\mathcal{T}$, hence by (\ref{mid1theo3.1}) we have
	\begin{equation}\label{midstep2_1}
	\int_{0}^{T}\int_{\Omega}f(\bar{m}_n)(m_{n+1}-\bar{m}_n)dxdt\le \int_{0}^{T}\int_{\Omega}f(\bar{m}_n)(m^{**}-\bar{m}_n)dxdt<-\frac{\ell}{2}.
	\end{equation}
	
    Recall that $\delta_n$ is the $n-$th learning rate in the fictitious play (\ref{fictitious}). Since $F\in C^{1,1}(\mathcal{C})$, $f$ is Lipschitz continuous, by (\ref{fictitious}) and Taylor expansions, for all $n\in\mathbb{N}$, we have
	\begin{equation}\label{midstep2_2}
	F(\bar{m}_{n+1})-F(\bar{m}_n)\le\int_{0}^{T}\int_{\Omega}\delta_nf(\bar{m}_n)(m_{n+1}-\bar{m}_n)dxdt+LC^2\delta_n^2,
	\end{equation}
	where $C:=2\left|\left|\tilde{m}\right|\right|_{\mathcal{C}}$ and $L$ is the Lipschitz constant of $f$. 
	
	Next we hope to show that, there exists a $\epsilon'< \epsilon$ which is small enough, such that when $n$ is large enough and $\|\bar{m}_n-m^*\|_{\mathcal{C}}<\epsilon'$, $F(\bar{m}_n)$ will decay to a value that is strictly less than $F(m^*)$, which contradicts with the assumption that $m^*$ is a cluster point. 
 
 By continuity of  $F(m)$ in the space $\mathcal{C}$, there exists $\epsilon'<\epsilon$, such that whenever $\left|\left|\bar{m}_n-m^*\right|\right|_{\mathcal{C}}<\epsilon'$, we have 
	\begin{equation}
        \label{midstep2_3}
	F(\bar{m}_n)\le F(m^*)+\frac{\ell(\epsilon-\epsilon')}{16C}.
	\end{equation}

 By the learning rate condition \eqref{cond_delta}, which we recall here for convenience
 $$
 \sum_n\delta_n=\infty,\quad\sum_n\delta_n^2<\infty,
 $$
 we can choose $N$ such that
 \begin{equation}\label{step2_bigN}
 \delta_n<\frac{(\epsilon-\epsilon')}{8C}\quad \text{for}\,\, n\ge N, \quad \text{and} \quad \sum_{n\ge N}\delta_n^2< \frac{\ell(\epsilon-\epsilon')}{16LC^3}.
 \end{equation}
 Since $m^*$ is a cluster point of $\bar{m}_n$ in the sense of $\|\cdot\|_{\mathcal{C}}$, we can find $N_0>N$ such that
$\left|\left|\bar{m}_{N_0}-m^*\right|\right|_{\mathcal{C}}<\epsilon'$. We define that
$$
M:=\min\{n>N_0:\sum_{k=N_0}^{n-1}\delta_k\ge\frac{7(\epsilon-\epsilon')}{8C}\}.
$$
By the fact that $\sum_n\delta_n=\infty$, $M$ is well-defined. We deduce from (\ref{step2_bigN}) that
$$
    \frac{7(\epsilon-\epsilon')}{8C}\le \sum_{k=N_0}^{M-1}\delta_k<\frac{\epsilon-\epsilon'}{C}.
$$

	Therefore, when $n=N_0+1,...,M$, $$\left|\left|\bar{m}_n-\bar{m}_{N_0}\right|\right|\le C\sum_{k=N_0}^{n-1}\delta_k<\epsilon$$ 
    holds. Hence by (\ref{midstep2_1}) , (\ref{midstep2_2}) and \ref{midstep2_3})  we have
	$$
	\begin{aligned}
	F(\bar{m}_M)-F(m^*)&\le 	F(\bar{m}_M)-F(\bar{m}_{N_0})+\frac{\ell(\epsilon-\epsilon')}{16C}\\&\le (-\frac{\ell}{2})\sum_{k=n}^{M-1}\delta_k+LC^2\sum_{k=n}^{M-1}\delta_k^2+\frac{\ell(\epsilon-\epsilon')}{16C}\\
	&< -\frac{7\ell(\epsilon-\epsilon')}{16C}+\frac{\ell(\epsilon-\epsilon')}{16C}+\frac{\ell(\epsilon-\epsilon')}{16C}\\
	&=-\frac{5\ell(\epsilon-\epsilon')}{16C}.
	\end{aligned}
	$$
	
	By step 1 and (\ref{midstep2_2}) we know that for all $n\in \mathbb{N}$, we have
	$$
		F(\bar{m}_{n+1})-F(\bar{m}_n)\le\int_{0}^{T}\int_{\Omega}\delta_nf(\bar{m}_n)(m_{n+1}-\bar{m}_n)dxdt+LC^2\delta_n^2\le LC^2\delta_n^2.
	$$
    In the last inequality, we applied \eqref{step1_result} with $m=\bar{m}_n$. Hence, for all $n\ge M$, we have the following inequality:
	$$
	F(\bar{m}_n)\le F(\bar{m}_M)+LC^2\sum_{n\ge M}\delta_n^2\le F(\bar{m}_M)+\frac{l(\epsilon-\epsilon')}{16C}\le F(m^*)-\frac{l(\epsilon-\epsilon')}{4C}.
	$$
	This is a contradiction with the continuity of $F$ and the assumption that $m^*$ is a cluster point. Thus the cluster point $m^*$ should satisfy
	\begin{equation}\label{midstep2_4}
	\left<\frac{\delta F}{\delta m}(m^*),m-m^*\right>\ge 0,
	\end{equation}
	for all $m\in\mathcal{T}$.
	\item We conclude that any cluster point $(u^*,m^*)$ is a mixed strategy equilibrium in this step. Suppose that  it is the subsequence $\{(u_{n_k},\bar{m}_{n_k})\}_{k=1}^\infty$ that converges to $(u^*,m^*)$.
 We first verify that $u^*$ satisfies:
 \begin{equation}\label{midstep3}
	\begin{cases}
	\max(-\partial _tu-\Delta u-f(m^*),u)=0,&(x,t)\in \Omega\times[0,T);\\
	u|_{t=T}=0,&x\in \Omega;\\
	u=0,&x\in \partial \Omega.
	\end{cases}
	\end{equation}
	
	Indeed, from the first equation of (\ref{fictitiousplay}), we know that ${u}_{n_k}$ is the solution to the following obstacle problem: 
	$$
	\begin{cases}
	\max(-\partial _tu-\Delta u-f(\bar{m}_{n_k-1}),u)=0,&(x,t)\in \Omega\times[0,T);\\
	u|_{t=T}=0,&x\in \Omega;\\
	u=0,&x\in \partial \Omega.
	\end{cases}
	$$
	Since $\delta_{n}\rightarrow0$, and $\bar{m}_{n_k}\to m^*$ by assumption, we know $\bar{m}_{n_k-1}\rightarrow m^*$ when $k\rightarrow\infty$. Thus by the continuity of $f$ with respect to $m$, as well as the uniqueness and stability of the obstacle problem (i.e. the map from $f$ to $u$ is also continuous, see section 5 in \cite{osmfgpde} for details),  $u^*$ is the solution of (\ref{midstep3}).

 It remains to verify
\begin{equation}\label{verifym1}
 \partial_t m^*-\Delta m^*=0\text{ in } \{u^*<0\},
 \end{equation}
 and
\begin{equation}\label{verifym2}
\int_{\{u^*=0\}}f(m^*)m^*\,dxdt=0.
\end{equation}
Define $m^{**}$ as in (\ref{mstarstar}). By (\ref{midstep2_4}) with $m=m^{**}$, and the relation \eqref{step1_result} with $\bar{m}_n:=m^*$, i.e.
$$
\left<\frac{\delta F}{\delta m}(m^*),m^{**}-m^*\right>\le 0,
$$
we know that
$$
\left<\frac{\delta F}{\delta m}(m^*),m^{**}-m^*\right>= 0.
$$
Then \eqref{verifym1} and \eqref{verifym2} follow from \eqref{cond1} and \eqref{cond2}, which are the conditions for the above equality to hold.
	
	\end{enumerate}
 Hence we have finished the proof.
\end{proof}
\begin{remark}
    The existence of the cluster point can be guaranteed through Lemma \ref{lem:weak_strong_comp}. To be precise, according to \eqref{unif_bound} in the proof of Theorem \ref{thm:regular_point}, each $m_{n}$ obtained by the fictitious play \eqref{fictitiousplay} is uniformly bounded in $\mathcal{B}$. Hence by Lemma \ref{lem:weak_strong_comp}, the sequence $\{\bar{m}_n\}_{n=1}^{\infty}$ admits at least a cluster point in $\mathcal{C}$, i.e. there exists a subsequence $\{\bar{m}_{n_k}\}_{k=1}^{\infty}\subset\{\bar{m}_n\}_{n=1}^{\infty}$, such that $\bar
    {m}_{n_k}$ converges to some $m^*$ in $\mathcal{C}$ as $k\rightarrow\infty$. Then following step 3 in the proof above, we can show that $(u_{n_k},\bar{m}_{n_k})$ converges to some $(u^*,m^*)$ in $\mathcal{A}\times \mathcal{C}$. 
\end{remark}
\subsection{Algorithm Based on Fictitious Play}
Theorem \ref{fictitious_convergence} provides the convergence result of the fictitious play as in Definition \ref{fictitiousplay}. We can turn the fictitious play into the following algorithm for finding the mixed strategy equilibrium.
\begin{algorithm}[htb]
		\caption{Fictitious Play to Find Equilibrium}
		\label{alg_fictitious}
		Given $u_0,m_0,m^0,f,\psi,itermax$. 
		
		$n:=0,\bar{m}_0:=m_0$.
		
		\textbf{While } ($n<itermax$)
				
		\quad Find $u_{n+1}$ by solving the obstacle problem (\ref{equation_for_u}) with $\mu:=\bar{m}_n$ and $u(x,T)=\psi(x,T,\bar{m}_n)$.
		
		\quad Find $m_{n+1}$ by solving the FP equation (\ref{equation_for_m}) with $u:=u_{n+1}$ and $m(x,0)=m^0(x)$.
		
		\quad $\bar{m}_{n+1}:=\delta_{n+1}m_{n+1}+(1-\delta_{n+1})\bar{m}_{n}$.
		
	    \quad$n:=n+1$.
		
		\textbf{End}
\end{algorithm}

In the remainder of this section, we will introduce the discretization method for the obstacle equation (\ref{equation_for_u}) and the Fokker-Planck equation (\ref{equation_for_m}) in algorithm \ref{alg_fictitious}. Assuming the spatial discretization grid size is $h$ and the time discretization step size is $\tau$, we use the notation $[u_n]_{j}^{k}$ to represent the numerical approximation of  $u_n(jh,k\tau)$, where $n$ denotes the iteration number in the fictitious play algorithm.

For the obstacle equation  (\ref{equation_for_u}), we need to solve it from $t=T:=M\tau$ backward to $t=0$. In order to use larger time steps, we consider the implicit scheme as follows:
\begin{equation}\label{implicit_obstacle}
		\max(\frac{[u_{n+1}]_j^k-[u_
  {n+1}]_j^{k+1}}{\tau}-\Delta_h[u_{n+1}]_j^k -[f(\bar{m}_n)]^k_j,[u_{n+1}]_j^k-[\psi(\bar{m}_n)]_j^k)=0,
\end{equation}
where $\Delta_h$ is the second-order central difference operator. We point out here that for each time step (\ref{implicit_obstacle}) is a discretized elliptic obstacle problem for $[u_{n+1}]^k$ since $\tau$ and $[u_{n+1}]^{k+1}$ are known. Numerical methods for this problem have been well studied. See \cite{obstacle_PLS},\cite{obstacle_MultigridFD},\cite{obstacle_SOR} for details. 

For the Fokker-Planck equation (\ref{equation_for_m}), we can also write down the implicit discretized scheme for it as follows:
\begin{equation}\label{implicit_FP}
\begin{cases}
\frac{[m_{n+1}]_j^{k}-[m_{n+1}]_j^{k-1}}{\tau}-\Delta_h[m_{n+1}]_j^{k}=0,&[u_{n+1}]_j^k<[\psi(\bar{m}_{n})]_j^k;\\
[m_{n+1}]_j^{k}=0,&[u_{n+1}]_j^k=[\psi(\bar{m}_{n})]_j^k. \\
\end{cases}
\end{equation}

However, iteration are unavoidable when numerically solving equations (\ref{implicit_obstacle}) and (\ref{implicit_FP}),  regardless of the method used. The application of nonlinear solvers renders the implicit scheme computationally inefficient. Therefore, in practice, semi-implicit schemes are preferred to reduce the computational cost. The semi-implicit scheme for $u$ can be written as follows:
\begin{equation}\label{semiimplicit_obstacle}
\begin{aligned}
	&\text{Step 1:}\frac{[\tilde{u}_{n+1}]_{j}^{k}-[u_{n+1}]_{j}^{k+1}}{\tau}-\Delta_h[\tilde{u}_{n+1}]_{j}^k -[f(\bar{m}_n)]^k_j=0;\\
	&\text{Step 2:}[u_{n+1}]_{j}^{k}=\min([\tilde{u}_{n+1}]_{j}^{k},[\psi(\bar{m}_n)]_j^k).
\end{aligned}
\end{equation}

And the semi-implicit scheme for $m$ can be written as:
\begin{equation}\label{semiimplicit_FP}
\begin{aligned}
&\text{Step 1:}\frac{[\tilde{m}_{n+1}]_j^{k}-[m_{n+1}]_j^{k-1}}{\tau}-\Delta_h[\tilde{m}_{n+1}]_j^{k}=0\\
&\text{Step 2:}[m_{n+1}]_j^{k}=[\tilde{m}_{n+1}]_j^{k}1_{[u_{n+1}]_j^k<[\psi(\bar{m}_n)]_j^k}.
\end{aligned}
\end{equation}

Semi-implicit schemes for $u$ and $m$ can be viewed as a two-step method that decouples the linear and nonlinear parts of each equation. In the first step, $u$ and $m$ are evolved using a standard implicit scheme on the whole domain. In the second step, a cutoff is applied to $u$ and $m$ respectively to account for the free boundary effects. We only need to solve a sparse system of linear equations for each time step in the semi-implicit schemes  (\ref{semiimplicit_obstacle}) and $(\ref{semiimplicit_FP}$.

Now we can summarize the finite difference algorithm into the following algorithm.
\begin{algorithm}[htb]
	\caption{Finite Difference Scheme for Mixed Strategy Equilibrium}
	\label{alg_FD}
	Given $u_0,m_0,m^0,f,\psi,itermax$. 
	
	$n:=0,\bar{m}_0:=m_0$.
	
	\textbf{While } ($n<itermax$)
	
	\quad Find $u_{n+1}$ by implicit scheme (\ref{implicit_obstacle}) (or semi-implicit scheme (\ref{semiimplicit_obstacle})).
	
	\quad Find $m_{n+1}$ by implicit scheme (\ref{implicit_FP}) (or semi-implicit scheme (\ref{semiimplicit_FP})).
	
	\quad $\bar{m}_{n+1}:=\delta_{n+1}m_{n+1}+(1-\delta_{n+1})\bar{m}_{n}$.
	
	\quad$n:=n+1$.
	
	\textbf{End}
\end{algorithm}

\subsection{Numerical Analysis}

In this part, our goal is to prove the convergence of algorithm \ref{alg_FD} when implicit scheme (\ref{implicit_obstacle}) and (\ref{implicit_FP}) are applied. The convergence analysis mirrors the proof for Theorem \ref{fictitious_convergence}, requiring only adapting the arguments to a discretized version. 

 For simplicity, we assume that the domain  $\Omega$ is the unit cubic $[0,1]^d$ in $\mathbb{R}^d$ and $\psi=0$. Discretize $\Omega$ with spatial scale $h=1/N$ and call the discretized domain $\Omega_h$. Denote $\Omega_h^{\circ}$ as the set of inner points of $\Omega_h$ and $\partial \Omega_h$ as the set of boundary points of $\Omega_h$.  Suppose the time step $\tau=T/K$. First, we introduce the implicit discretized system for mixed strategy equilibrium.
\begin{definition}(implicit discretized system for mixed strategy equilibrium)\label{implicit_mixed_strategy}
	We define $$\Gamma(u):=\{(j,k):\frac{u^{k}_{j}-u^{k+1}_{j}}{\tau}-\Delta_hu^{k}_{j} -f(m)^k_j<0\}.$$ A couple $(u,m)\in \mathbb{R}^{(K+1)(N+1)^d}\times\mathbb{R}^{(K+1)(N+1)^d}$ is a solution of the discretized system for the mixed strategy equilibrium if
	\begin{align}\label{implicit_mixed_strategy_pde}
\begin{cases}
\max(\frac{u^{k}_{j}-u^{k+1}_{j}}{\tau}-\Delta_hu^{k}_{j} -f(m)^k_j,u_{j}^k)=0, &\text{in }\{0,...,K-1\}\times\Omega_h^{\circ};\\
\frac{m_{j}^{k}-m_{j}^{k-1}}{\tau}-\Delta_hm_{j}^{k}=0,&\text{in }\{u^k_{j}<0\};\\
\frac{m_{j}^{k}-m_{j}^{k-1}}{\tau}-\Delta_hm_{j}^{k}\le0,\,m_j^k\ge 0&\text{in }\{1,...,K\}\times\Omega_h^{\circ};\\
m^0_j=[m^0]_j,\,u^K_j=0,&\text{in }\Omega_h;\\
u^k_j=0,&\text{in }\{0,...,K-1\}\times\partial \Omega_h;\\
m^k_j=0,&\text{in }\{1,...,K\}\times\partial \Omega_h;\\
\sum_{\Gamma(u)}(f(m)_j^km_j^k)=0.
\end{cases}
\end{align}
\end{definition}
\begin{remark}
We note here that $\Gamma(u)\subset \{u_{j}^k=0\}$. When there exists some point such that,
\begin{equation}\label{singular}
\begin{cases}
\frac{u^{k}_{j}-u^{k+1}_{j}}{\tau}-\Delta_hu^{k}_{j} -f(m)^k_j=0;\\
u_j^k=0;\\
f(m)^k_j>0,
\end{cases}
\end{equation}
the complementary condition 
\begin{equation}\label{dis_comple}
\sum_{\Gamma(u)}(f(m)_j^km_j^k)=0
\end{equation}
will be weaker than the following one
\begin{equation}\label{cont_comple}
\sum_{\{u_j^k=0\}}(f(m)_j^km_j^k)=0.
\end{equation}
However, a point satisfies (\ref{singular}) should lie on $\partial\{u_j^k=0\}$ and thus, the condition (\ref{dis_comple}) is effectively equivalent to (\ref{cont_comple}) when the free boundary is regular. 
\end{remark}

Before stating the main result, we present a property of the implicit scheme for obstacle equations: the discretized solution $u$ continuously depends on the discretized source term $f$.

\begin{lemma}\label{cont_dep_obs}
Consider the following discrete obstacle problem:
\begin{equation}\label{discrete_obstacle}
\begin{cases}
\max(\frac{u^{k}_{j}-u^{k+1}_{j}}{\tau}-\Delta_hu^{k}_{j} -f^k_j,u_{j}^k)=0, &\text{in }\{0,...,K-1\}\times\Omega_h^{\circ};\\
u^K_j=0,&\text{in }\Omega_h;\\
u^k_j=0,&\text{in }\{0,...,K-1\}\times\partial \Omega_h.
\end{cases}
\end{equation} 
Suppose that $f,\,f'\in\mathbb{R}^{K(N-1)^d}$. For any $\epsilon>0$, when $\|f-f'\|_{\infty}<\epsilon$, the difference between two corresponding solutions satisfies $\|u-u'\|_{\infty}<\epsilon T$.
\end{lemma}
\begin{proof}
Without loss of generality, we assume $d=1$ and denote $M=\#\Omega_h^0$, $[u]^k=(u_1^{k},...,u_M^{k})^t,\,[u']^k=([u']_1^{k},...,[u']_M^{k})^t,\,[f]^k=(f_1^{k},...,f_M^{k})^t,\,[f']^k=([f']_1^{k},...,[f']_M^{k})^t$, for $k=0,\,\cdots,\,K-1,$ where $K\tau=T$. 

For $k\in \{0,\,\cdots,\,K-1\}$, we suppose that $\ell=\arg\max_j|[u]_j^k-[u']_j^k|$. Without loss of generality, we assume that $[u]_\ell^k>[u']_{\ell}^k$. We first show that
 \begin{equation}\label{eq:cont_uupidiff_final}
    \|[u]^k-[u']^k\|_{\infty}\le \|[u]^{k+1}-[u']^{k+1}\|_{\infty}+\tau \|f-f'\|_{\infty}.
\end{equation}

\begin{itemize}
    \item If $0>[u]_\ell^k>[u']_\ell^k$, then according to the obstacle equation \eqref{discrete_obstacle}, we have
    \begin{equation}\label{eq:cont_ueq}
        {[u]}_\ell^k-[u]_\ell^{k+1}-\tau\Delta_h[u]_{\ell}^k-\tau[f]_{\ell}^k=0,
    \end{equation}
    and
    \begin{equation}\label{eq:cont_upieq}
        {[u']}_\ell^k-[u']_\ell^{k+1}-\tau\Delta_h[u']_{\ell}^k-\tau[f]_{\ell}^k=0.
    \end{equation}
    We subtract \eqref{eq:cont_ueq} from \eqref{eq:cont_upieq} and obtain
    \begin{equation}\label{eq:cont_uupidiff}
        {[u-u']}_\ell^k-[u-u']_\ell^{k+1}-\tau\Delta_h[u-u']_{\ell}^k-\tau[f-f']_{\ell}^k=0.
    \end{equation}
    Since for any $i\in\{0,\,\cdots,\,M+1\}$, We have $[u-u']_\ell^k=|[u]_\ell^k-[u']_\ell^k|=\max_j|[u]_j^k-[u']_j^k|\ge [u-u']_i^k$. Therefore,
    \begin{equation}\label{eq:cont_lapuupi}
       -\tau\Delta_h[u-u']_{\ell}^k=\frac{\tau}{h^2}\left(2[u-u']_{\ell}^k-[u-u']_{\ell+1}^k-[u-u']_{\ell-1}^k\right)\ge 0. 
    \end{equation}
    We utilize the equality \eqref{eq:cont_lapuupi} in \eqref{eq:cont_uupidiff} and conclude that
    \begin{equation}\label{eq:cont_uupidiff_simp}
    {[u-u']}_\ell^k-[u-u']_\ell^{k+1}-\tau[f-f']_{\ell}^k\le 0.
    \end{equation}
    By the fact that ${[u-u']}_\ell^k=\|[u]^k-[u']^k\|_{\infty}$,\,${[u-u']}_\ell^{k+1}\le \|[u]^{k+1}-[u']^{k+1}\|_{\infty}$,\,and ${[f-f']}_\ell^k=\|[f]^k-[f']^k\|_{\infty}\le \|f-f'\|_{\infty}$, inequality \eqref{eq:cont_uupidiff_simp} can be further estimated as
    \eqref{eq:cont_uupidiff_final}.

    \item $0=[u]_\ell^k>[u']_\ell^k$, the obstacle equation \eqref{eq:cont_ueq} for $[u]_{\ell}^k$ is replaced to
    \begin{equation}
         {[u]}_\ell^k-[u]_\ell^{k+1}-\tau\Delta_h[u]_{\ell}^k-\tau[f]_{\ell}^k\le 0,
    \end{equation}
    and we still have the obstacle equation \eqref{eq:cont_upieq} for $[u']_{\ell}^k$. Then proceeding the remainder of the proof for the case $0>[u]_{\ell}^k>[u']_{\ell}^k$, we still have the estimation \eqref{eq:cont_uupidiff_final}.
\end{itemize}
 For any $k\in\{0,\,\cdots,\,K-1\}$, we add up inequality \eqref{eq:cont_uupidiff_final} for $k,\,k+1,\,\cdots,K-1$ and obtain
    $$
       \|[u]^k-[u']^k\|_{\infty}\le   \|[u]^K-[u']^K\|_{\infty}+(K-k)\tau\|f-f'\|_{\infty}=(K-k)\tau\|f-f'\|_{\infty}\le T\epsilon,
    $$
    yielding that
    \begin{equation}\label{eq:cont_uupidiff_addup}
    \|u-u'\|_{\infty}\le T\epsilon.
    \end{equation}
\end{proof}

Now we can state the main convergence result in this section.
\begin{theorem}\label{convergence_implicit}
		Given $\tau,\,h>0$. Consider the Euclidean space $\mathbb{R}^{K(N-1)^d}$ equipped with the norm $\|\cdot\|$ such that
  $$
  \|g\|:=\sqrt{ \tau h^d\sum_{k=0}^{K-1}\sum_{j\in\Omega_h}(g_j^k)^2}
  $$
  for all $g\in \mathbb{R}^{K(N-1)^d}$. Assume that $m\in\mathbb{R}^{K(N-1)^d}$, $f(m):\mathbb{R}^{K(N-1)^d}\rightarrow\mathbb{R}^{K(N-1)^d}$ is a Lipschitz continuous vector function, and there exists a function $\Phi(m)\in C^1(\mathbb{R}^{K(N-1)^d};\mathbb{R})$ such that 
  $$f(m)=\frac{\partial\Phi}{\partial m},\text{ i.e. }[f(m)]_j^k=\frac{\partial\Phi}{\partial m_{j}^{k}}$$
  for any $j\in\{1,...,N-1\},k\in\{0,...,K-1\}$. Then any cluster point of the sequences $(u_n,\bar{m}_n)$ obtained by algorithm \ref{alg_FD} with the implicit scheme would be a solution of ($\ref{implicit_mixed_strategy_pde}$). 
\end{theorem}
\begin{proof}
	The spirit of the proof is analog to the one in theorem \ref{fictitious_convergence}. We divide the proof into 3 steps just parallel to the proof of theorem \ref{fictitious_convergence}.
	\begin{enumerate}
		\item 
        Define the test set $\mathcal{T}_{h,\tau}$  as following:
        $$
        \mathcal{T}_{h,\tau}:=\{m\in\mathbb{R}^{K(N+1)^d}:m_j^k=0\text{ when }j\in \partial \Omega_h,\, m_j^k\ge 0,\,\frac{m_{j}^{k}-m_{j}^{k-1}}{\tau}-\Delta_hm_{j}^{k}\le0,\,m^0_j=[m^0]_j\}.
        $$
		We denote $\left<\cdot,\cdot\right>$ as the inner product in $\mathbb{R}^{K(N-1)^d}$ with respect to the norm $\|\cdot\|$. In this step, we aim to show that for any  $m\in\mathcal{T}_{h,\tau}$, we have 
		\begin{equation}\label{discrete_variation}
				\left<\frac{\partial \Phi}{\partial m},m_{n+1}-m\right>=\tau h^d\sum_{k=0}^{K-1}\sum_{j\in\Omega_h}f(\bar{m}_n)_j^k[m_{n+1}-m]_j^k\le 0.
		\end{equation} 
		 Indeed, according to algorithm \ref{alg_FD}, in $\Gamma(u_{n+1})$,
		 we have $[m_{n+1}]^k_{j}=0$, and thus $[m_{n+1}-m]_j^k\le 0$.
		 Therefore we have 
        $$
		 [\frac{[u_{n+1}]^{k}_{j}-[u_{n+1}]^{k+1}_{j}}{\tau}-\Delta_h[u_{n+1}]^{k}_{j}-[f(\bar{m}_n)]^k_j][m_{n+1}-m]_j^k\ge 0.
        $$
        
        In $\Gamma(u_{n+1})^c$, we have$$
    \frac{[u_{n+1}]^{k}_{j}-[u_{n+1}]^{k+1}_{j}}{\tau}-\Delta_h[u_{n+1}]^{k}_{j} -[f(\bar{m}_n)]^k_j=0,
        $$
        and thus
        $$
        [\frac{[u_{n+1}]^{k}_{j}-[u_{n+1}]^{k+1}_{j}}{\tau}-\Delta_h[u_{n+1}]^{k}_{j} -[f(\bar{m}_n)]^k_j][m_{n+1}-m]_j^k= 0.
        $$
		Hence,
		\begin{equation}\label{discrete_mid1}
		\begin{aligned}
		&\sum_{k=0}^{K-1}\sum_{j\in\Omega_h}[f(\bar{m}_n)]_j^k[m_{n+1}-m]_j^k\\
		\le &\sum_{k=0}^{K-1}\sum_{j\in\Omega_h}[\frac{[u_{n+1}]^{k}_{j}-[u_{n+1}]^{k+1}_{j}}{\tau}-\Delta_h[u_{n+1}]^{k}_{j} ][m_{n+1}-m]_j^k\\
		= &\sum_{k=1}^{K}\sum_{j\in\Omega_h}[u_{n+1}]_{j}^k[\frac{[m_{n+1}-m]_j^k-[m_{n+1}-m]_j^{k-1}}{\tau}-\Delta_h[m_{n+1}-m]_j^{k}]\\
    = &\sum_{(j,k)\in\{u_j^k<0\}}[u_{n+1}]_{j}^k[\frac{[m_{n+1}-m]_j^k-[m_{n+1}-m]_j^{k-1}}{\tau}-\Delta_h[m_{n+1}-m]_j^{k}]\\
		\end{aligned}
		\end{equation}
  with equality if and only if 
  $$
  m_j^k=[m_{n+1}]_j^k,\text{ in }\Gamma(u_{n+1}),
  $$
  which implies (recall $[m_{n+1}]_j^k=0$ in $\Gamma(u_{n+1})$)
  \begin{equation}\label{cond1dis}
  \sum_{(j,k)\in \Gamma(u_{n+1})}[f(\bar{m}_n)]_j^km_j^k=0.
  \end{equation}
  To prove (\ref{discrete_variation}), it suffices to show that the right hand side of (\ref{discrete_mid1}) is not greater than 0. We know from algorithm \ref{alg_FD} that when $u^k_{j,n+1}<0$, 
		$$\frac{[m_{n+1}]_{j}^k-[m_{n+1}]_{j}^{k-1}}{\tau}-\Delta_h[m_{n+1}]_{j}^k=0$$
		holds. 
  Hence for all $m\in\mathcal{T}_{h,\tau}$, we have $$
		\frac{[m_{n+1}-m]_j^k-[m_{n+1}-m]_j^{k-1}}{\tau}-\Delta_h[m_{n+1}-m]_j^{k}\ge0.$$
		And thus  $$
		[u_{n+1}]_{j}^k[\frac{[m_{n+1}-m]_j^k-[m_{n+1}-m]_j^{k-1}}{\tau}-\Delta_h[m_{n+1}-m]_j^{k}]\le 0,
		$$
  which yields that the right hand side of (\ref{discrete_mid1}) is not greater than 0, with equality if and only if 
  \begin{equation}\label{cond2dis}
  \frac{m_{j}^{k}-m_{j}^{k-1}}{\tau}-\Delta_hm_{j}^{k}=0,\text{ in }\{[u_{n+1}]^{k}_{j}<0\};
  \end{equation}
		
		\item 	Suppose that  $(u^*,m^*)$ is a cluster point of the sequence $(u_n,\bar{m}_n)$ obtained by algorithm \ref{alg_FD}, following the same argument in step 2 in theorem  \ref{fictitious_convergence}, we know that $m^*$ should be the minimizer of
		$$
		\min_{
			m\in\mathcal{T}_{h,\tau}}
		\tau h^d\sum_{k=0}^{K-1}\sum_{j\in\Omega_h}[f(m^*)]_j^km_j^k.
		$$

          \item We conclude that any cluster point $(u^*,m^*)$ is a solution to (\ref{implicit_mixed_strategy_pde}).
         We first verify that $u^*$ will satisfy the discretized obstacle problem as follows:
          \begin{equation}\label{step3obstacle}
        \begin{cases}
        \max(\frac{u^{k}_{j}-u^{k+1}_{j}}{\tau}-\Delta_hu^{k}_{j} -[f(m^*)]_j^k,u_{j}^k)=0, &\text{in }\{0,...,K-1\}\times\Omega_h^{\circ};\\
        u^K_j=0,&\text{in }\Omega_h;\\
        u^k_j=0,&\text{in }\{0,...,K-1\}\times\partial \Omega_h.
        \end{cases}
        \end{equation}
        Indeed, we know that $\bar{u}_{n_k}$ is the solution of the following discretized obstacle problem:
        $$
        \begin{cases}
        \max(\frac{u^{k}_{j}-u^{k+1}_{j}}{\tau}-\Delta_hu^{k}_{j} -[f(\bar{m}_{n_k-1})]_j^k,u_{j}^k)=0, &\text{in }\{0,...,K-1\}\times\Omega_h^{\circ};\\
        u^K_j=0,&\text{in }\Omega_h;\\
        u^k_j=0,&\text{in }\{0,...,K-1\}\times\partial \Omega_h.
        \end{cases}
        $$
        Since $\delta_n\rightarrow0$, we know $\bar{m}_{n_k-1}\rightarrow m^*$ when $k\rightarrow\infty$. Thus by the continuity of $f$ with respect to $m$ and lemma \ref{cont_dep_obs}, $u^*$ is the solution of (\ref{step3obstacle}).
        The equality satisfied by $m^*$ follows from conditions (\ref{cond1dis}) and (\ref{cond2dis}), using the same arguments as step 3 in the proof of theorem \ref{fictitious_convergence}. For brevity, we omit the details.
\end{enumerate}
\end{proof}

\begin{remark}
A similar convergence result holds for the explicit scheme with the CFL condition.
\end{remark}

\section{Numerical Experiments}
In this section, we conduct several numerical experiments to demonstrate the effectiveness of the proposed semi-implicit finite difference algorithm (Algorithm \ref{alg_FD}). Through these experiments, we examine the convergence properties of our algorithm highlighting the implementation of the fictitious play. We demonstrate that our requirement for $\delta_n$ in \eqref{cond_delta} is a sufficient yet unnecessary condition. However, in certain cases where the pure strategy equilibrium may not exist, the iteration method may not converge if the condition \eqref{cond_delta} is violated. 

\subsection{A Non-local OSMFG Example}
\label{sec:nonlocal}

\paragraph{Setup.}
In this example, the state of the representative agent belongs to the domain $[0, 1] \times \mathbb{R}$. It dynamic is a Brownian motion, i.e.
\begin{equation*}
    X_t=X_0+\sqrt{2}W_t
\end{equation*}

The initial population distribution $m^0$ is in a Gaussian form: 
\begin{equation*}
    m^{0}(x)=\frac{1}{2\sqrt{2\pi}}\exp{\left(-\frac{x^2}{8}\right)}
\end{equation*}

The running cost $f(x,t,m)$ is defined as
\begin{equation*}
    f(x,t,m)=x-\frac{\int_\Omega \xi m(\xi,t) d\xi}{\int_\Omega  m(\xi,t) d\xi}
\end{equation*}
and the stopping cost is defined as
\begin{equation*}
    \psi(x,t,m)=-t
\end{equation*}

Intuitively, as the agent aims to minimize expected cost, the running cost encourages the agent to remain in the game when its state is below the average state of other remaining agents, while the stopping cost encourages the agent to continue playing for a longer duration.

It is easy to check that the cost functions above are equivalent with
\begin{equation*}
        \tilde{f}(x,t,m)=x-1-\frac{\int_\Omega \xi m(\xi,t) d\xi}{\int_\Omega  m(\xi,t) d\xi},\quad \tilde{\psi}(x,t,m)=0
\end{equation*}

With $\tilde{f}$ and $\tilde{\psi}$, we can formulate this problem into the PDE form of Equation (\ref{fictitious}) and discretize it as Equations (\ref{semiimplicit_obstacle}) and (\ref{semiimplicit_FP}).
In this experiment, we approximate the unbounded domain $\Omega=\mathbb{R}$ by a bounded domain $\hat{\Omega}=[-5,5]$ with Dirichlet boundary conditions and discretize it uniformly with the mesh $\Omega_{h}=\{-5,-5+h,\dots,5\}$. The time step is set to $\tau= h^2$ for the sake of accuracy. 
We remark that with the semi-implicit scheme, numerical experiments show that much larger time steps can be taken to produce convergent results. We omit to report the standard convergence tests with respect to the discretization error of the PDEs, and the following tests are mainly devoted to exploring the iteration scheme in the fictitious play. 

\paragraph{Numerical Result.}

\begin{figure}
\centering
\subfigure[Distribution $m$]{
    \centering
    \includegraphics[width=0.48\linewidth]{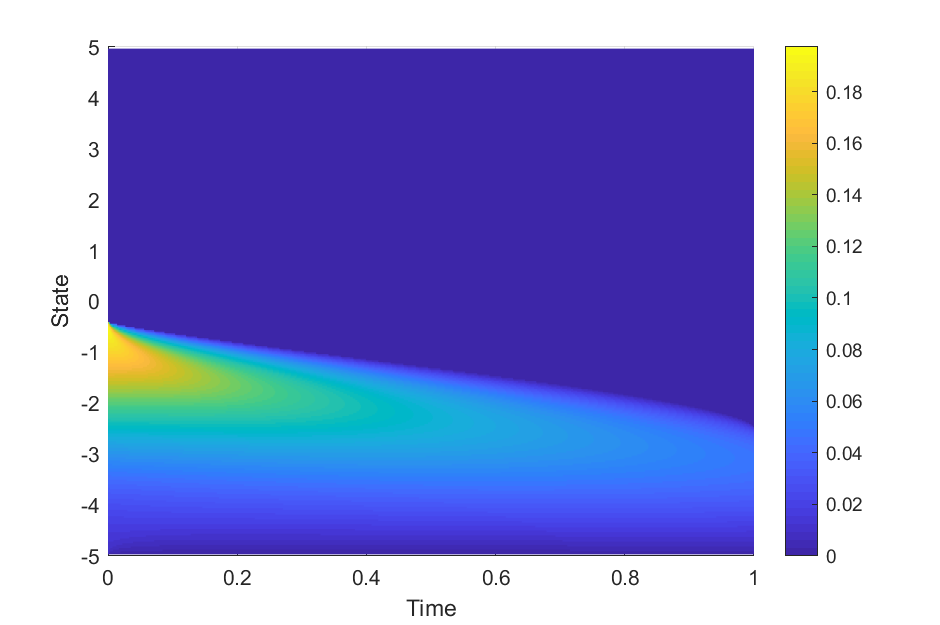}
    \label{fig1a}
}
\subfigure[Agents stay in the game]{
    \centering
    \includegraphics[width=0.48\linewidth]{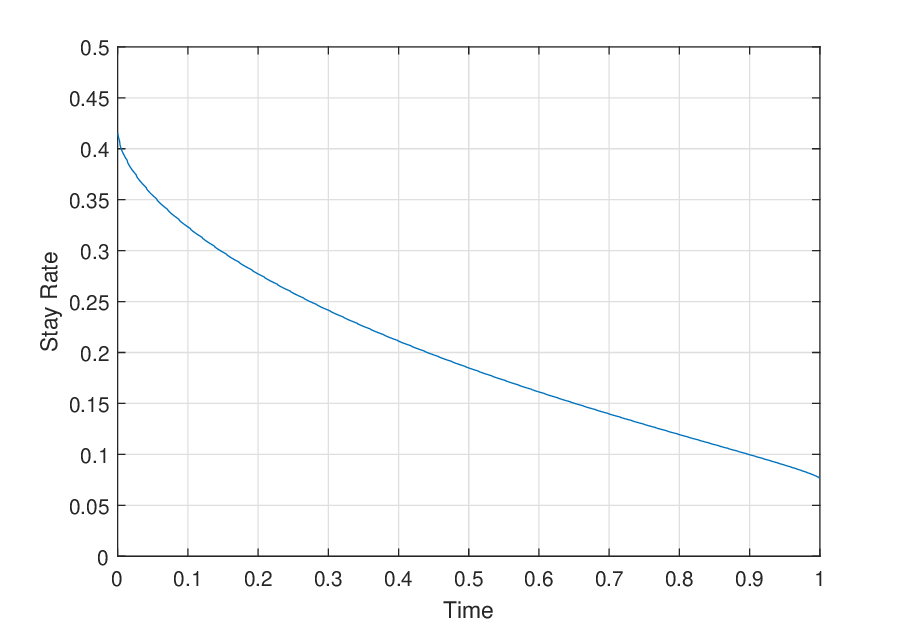}
    \label{fig1c}
}\\
\subfigure[Value function $u$]{
    \centering
    \includegraphics[width=0.48\linewidth]{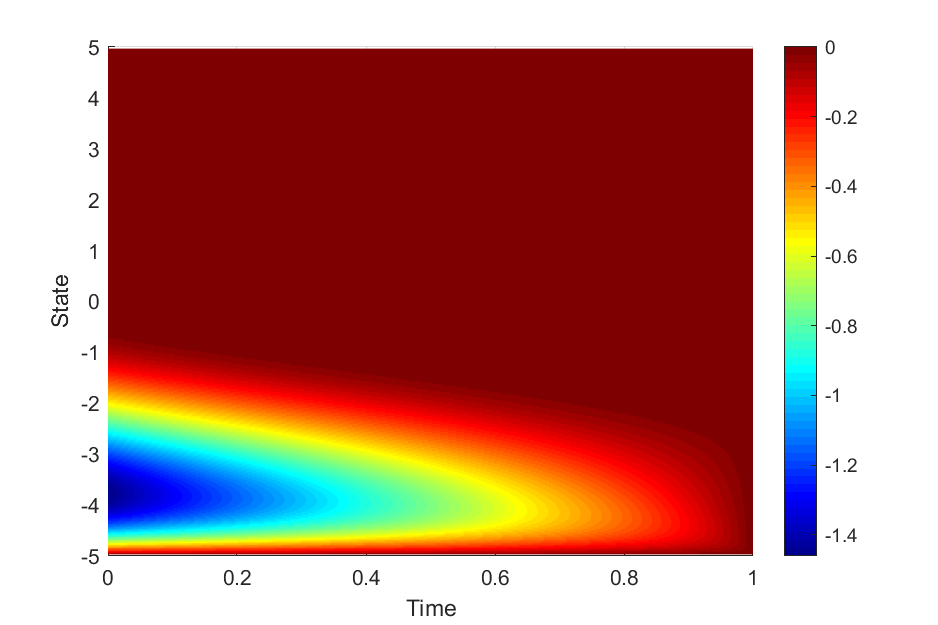}
    \label{fig1b}
}
\subfigure[Exiting boundary]{
    \centering
    \includegraphics[width=0.48\linewidth]{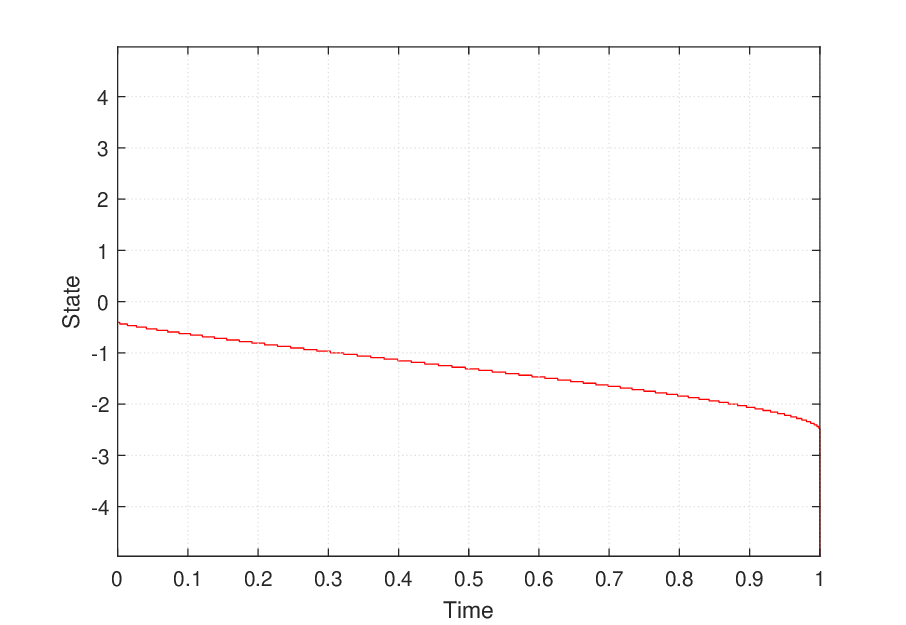}
    \label{fig1d}
}
\caption{
Numerical results for the example given in Section 4.1, $h=2^{-5}$, $\delta_n=1/n$, 1000 iterations.
\textbf{Top-Left} shows the evolution of the distribution $m$. 
\textbf{Top-Right} shows the rate of agents remaining in the game as time progresses. 
\textbf{Bottom-Left} shows the evolution of the value function $u$.
\textbf{Bottom-Right} highlights the exiting boundary.
}
\label{fig1}
\end{figure}

Figure \ref{fig1} demonstrates a numerical solution of this example. Here, the mesh size is set to $h=2^{-5}$ and the learning rate is $\delta_n=1/n$. The solution is computed for 1000 iterations to obtain the result shown. As depicted in figure \ref{fig1a}, agents starting with a high state exit the game immediately, while agents starting with a low state remain in the game for a longer period of time. Figure \ref{fig1c} shows the amount of agents remaining in the game as time progresses. Figure \ref{fig1b} illustrates the value function $u$ and Figure \ref{fig1d} illustrates the exiting boundary, i.e. the boundary of $\{u=0\}$. Once active agents reach this boundary for the first time, they exit the game immediately. 

Since an analytical solution could not be obtained, here we use two metrics to evaluate the convergence of our algorithm. 
The first metric utilized the numerical result on finer grids as an approximation of the true solution. Specifically, the numerical result on a grid of $h=2^{-7}$ with $\delta_n=1$ after 1000 iterations was used as the baseline solution. The error was defined as the $l$-2 norm between the numerical solution $\bar{m}_n$ and the baseline solution $\bar{m}$:
\begin{equation*}
\epsilon_n= \left(\tau h \sum_{k=0}^M\sum_{j\in\Omega_h}(\bar{m}^k_{j,n}-\bar{m}^k_{j})^2\right)^{1/2}
\end{equation*}

The second metric for convergence is the exploitability of the solution. Exploitability is a concept introduced in \cite{Perrin2020FictitiousPF}. It quantifies the average gain of the representative agent by switching to the optimal policy while the other agents retain their original policies. In our case, the exploitability can be defined as
\begin{equation*}
e_n=\tau h\sum_{k=0}^M\sum_{j\in\Omega_h}\tilde{f}(x_j,t_k,\bar{m}_n)(m^k_{j,n+1}-\bar{m}^k_{j,n})
\end{equation*}
where $m_{n+1}$ is the updated distribution, as defined in Algorithm \ref{alg_FD}.

\begin{figure}
\centering
\subfigure[$l$-2 error]{
    \centering
    \includegraphics[width=0.48\linewidth]{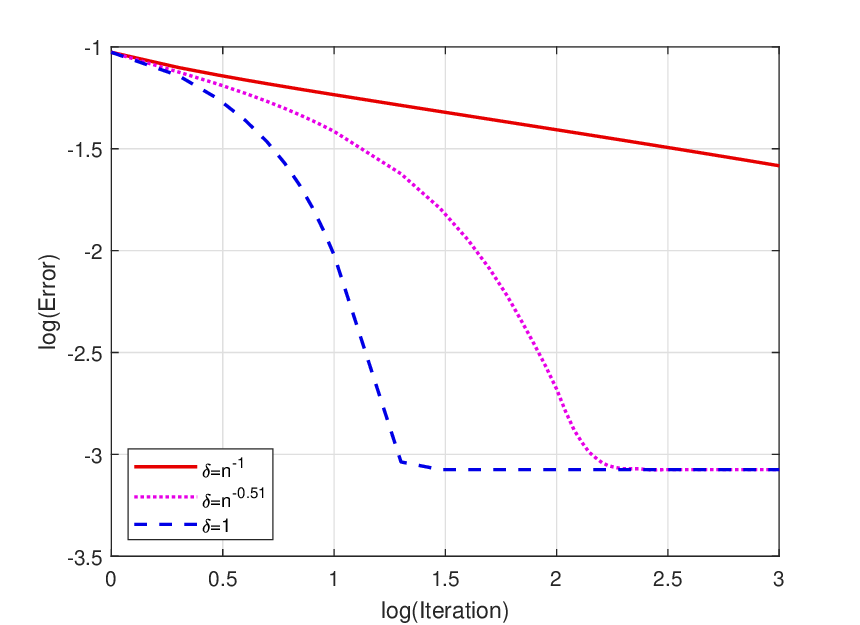}
    \label{fig3a}
}
\subfigure[Exploitability]{
    \centering
    \includegraphics[width=0.48\linewidth]{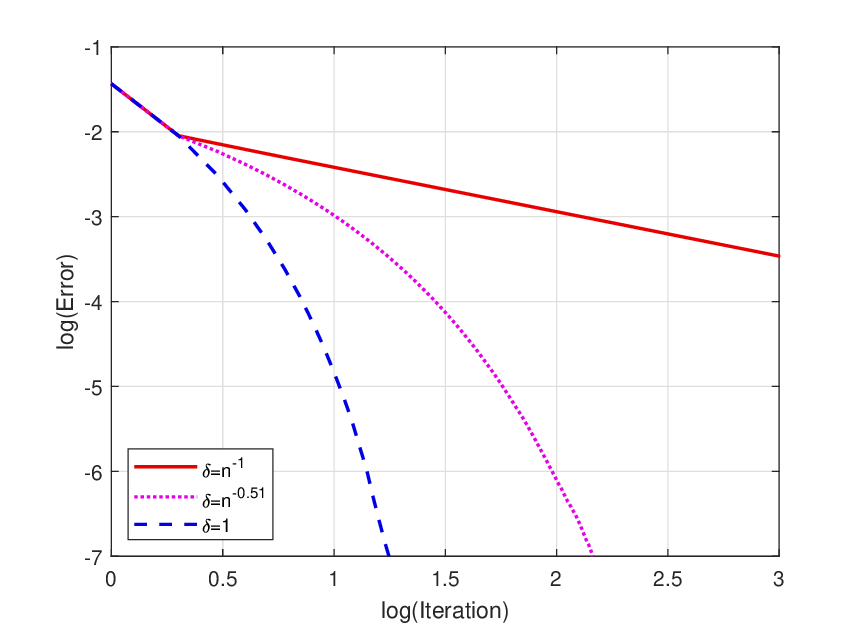}
    \label{fig3b}
}
\caption{
\centering
$\log_{10}$-$\log_{10}$ plots, error vs. iteration with different learning rate $\delta_n$. 
\textbf{Left} shows the $l$-2 errors with respect to the baseline solution. 
\textbf{Right} shows the explotability.
}
\label{fig3}
\end{figure}

In Figure \ref{fig3}, log-log plots of error versus iteration numbers are shown for different learning rates $\delta_n$. For this example, taking $\delta_n=1$ leads to the fastest convergence, although with $\delta_n=1$ no fictitious play is implemented. The reason for this is that the convergence requirement for the learning rate $\delta_n$ in Definition \ref{fictitiousplay} is only a sufficient condition, rather than a necessary one. For problems with a pure strategy equilibrium, the condition in definition \ref{fictitiousplay} can be relaxed and a more aggressive learning rate can be taken to obtain faster convergence. However, we would like to stress that we could not know a priori whether a pure strategy equilibrium exists, but the proposed algorithm, which is based on fictitious play, always produces convergence results.

\begin{figure}
\centering
\subfigure[$l$-2 error]{
    \centering
    \includegraphics[width=0.48\linewidth]{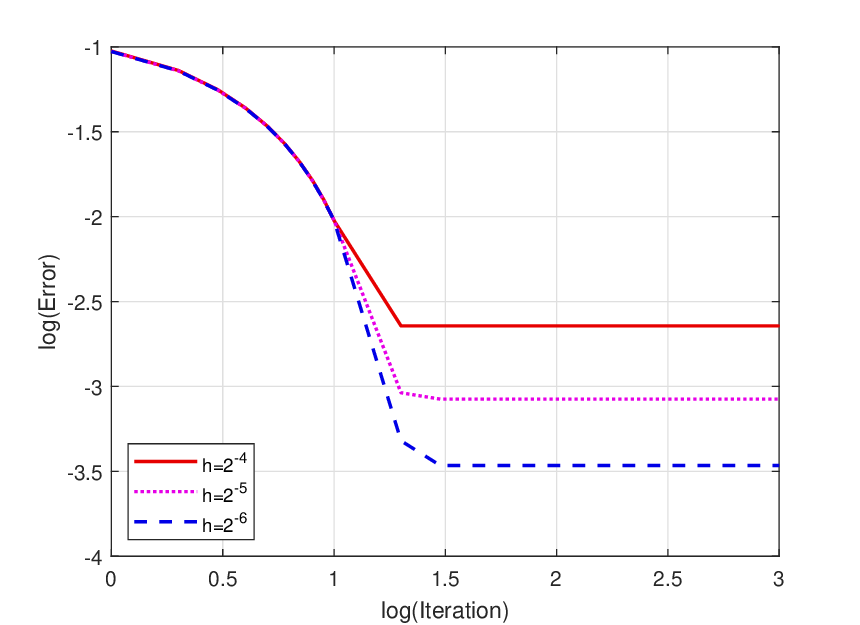}
    \label{fig2a}
}
\subfigure[Exploitability]{
    \centering
    \includegraphics[width=0.48\linewidth]{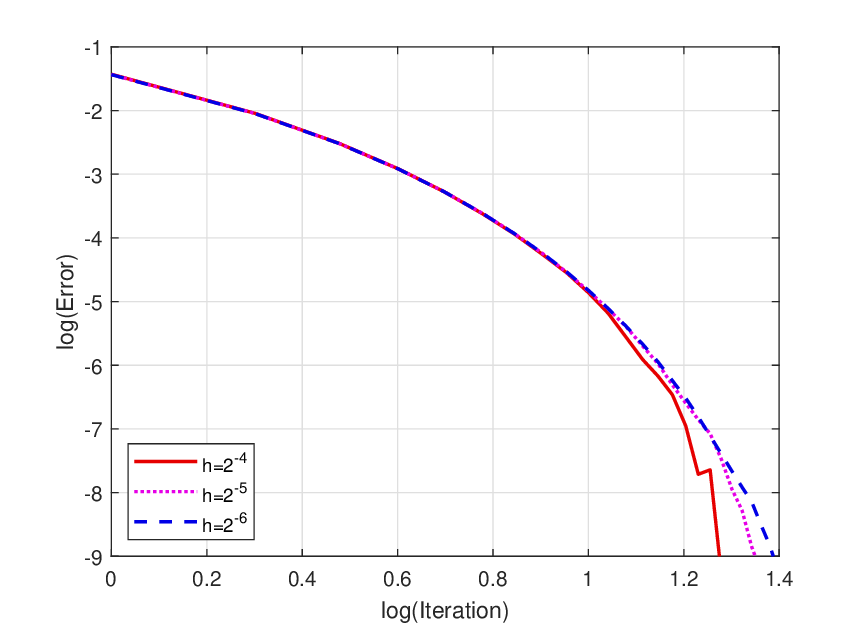}
    \label{fig2b}
}
\caption{\centering
$\log_{10}$-$\log_{10}$ plots, error vs. iteration with different mesh size $h$. 
\textbf{Left} shows the $l$-2 errors with respect to the baseline solution. 
\textbf{Right} shows the explotability.
}
\label{fig2}
\end{figure}

In Figure \ref{fig2}, we show log-log plots of error versus iteration number for different mesh sizes $h$. For these plots, the learning rate is set to $\delta_n=1$. From the plots, we can observe that the exploitabilities are very close across the different mesh sizes $h$. This means that exploitability is not affected by the mesh size among the test set, thus the discretization error is not a dominating factor in this test. On the other hand, the $l$-2 errors reach plateaus after about 20 iterations, which we interpret as the numerical solution has reached the mesh resolution limit respectively, and more iterations can no longer help to reduce the overall numerical error. In this case, the converged $l$-2 error scales as approximately $O(h^{1.37})$.

\subsection{A 2D Local OSMFG Example}
\label{sec:2d_local}
 
\paragraph{Setup.}
This example is a two-dimensional analogue of a local congestion OSMFG: the running cost depends on $m$ only through its local value, and the stopping cost is independent of $m$. The setting therefore falls into the potential-game framework of Sections~2--3. The state of the representative agent is $X_t=(X_t^1,X_t^2)\in\Omega=[-5,5]^2$, with terminal time $T=0.6$. The law of motion is the two-dimensional Ornstein--Uhlenbeck process
\begin{equation*}
    dX_t = -X_t\,dt + \sqrt{2}\,dW_t,
\end{equation*}
where $W_t$ is a standard two-dimensional Brownian motion.

The initial distribution is a bimodal Gaussian, normalized to total mass $2$:
\begin{equation*}
    m^{0}(x_1,x_2)=C\left[\exp\!\left(-\frac{(x_1-1.5)^2+x_2^2}{2}\right)
    +\exp\!\left(-\frac{(x_1+1.5)^2+x_2^2}{2}\right)\right],
\end{equation*}
where $C>0$ is chosen so that $\int_\Omega m^{0}=2$.
The running cost is local in $m$,
\begin{equation*}
    f(x,t,m)=-0.2+3\,m(x,t),
\end{equation*}
and the stopping cost is constant,
\begin{equation*}
    \psi(x,t,m)=0.
\end{equation*}
The drift is discretized by a first-order upwind scheme. The computational domain $[-5,5]^2$ is equipped with Dirichlet boundary conditions and discretized with mesh size $h=5\cdot 2^{-4}$ and time step $\tau=0.12\,h^2$. Each fictitious-play iteration consists of one backward obstacle solve and one forward Fokker--Planck solve, as in Algorithm~\ref{alg_FD}, with no additional modification of the value function.

The inward drift concentrates the two initial peaks toward the origin. The congestion term $3m$ penalizes overcrowding and encourages exit, while a drop in density makes remaining in the game attractive again. This anti-monotone coupling is not expected to admit a pure strategy equilibrium.

\paragraph{Numerical Result.}

\begin{figure}
\centering
\subfigure[$m(x,t{=}0)$]{
    \centering
    \includegraphics[width=0.48\linewidth]{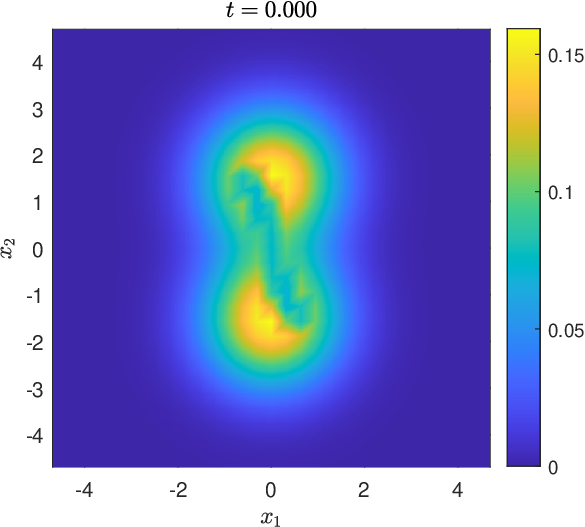}
    \label{fig4a}
}
\subfigure[$m(x,t{=}T/3)$]{
    \centering
    \includegraphics[width=0.48\linewidth]{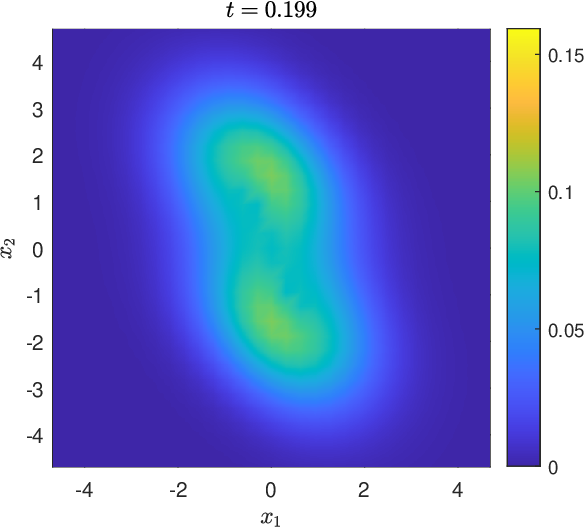}
    \label{fig4b}
}\\
\subfigure[$m(x,t{=}2T/3)$]{
    \centering
    \includegraphics[width=0.48\linewidth]{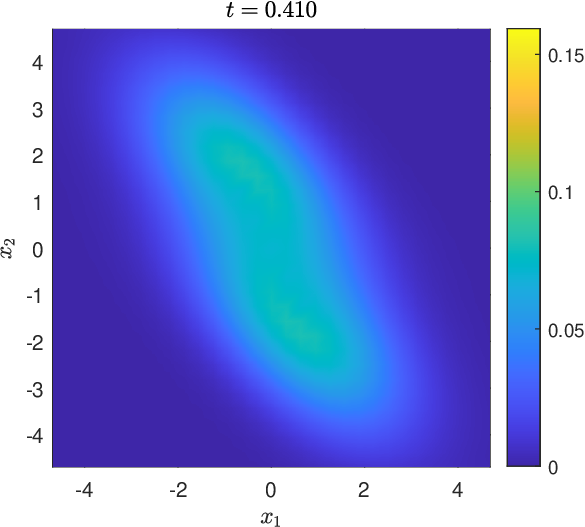}
    \label{fig4c}
}
\subfigure[$m(x,t{=}T)$]{
    \centering
    \includegraphics[width=0.48\linewidth]{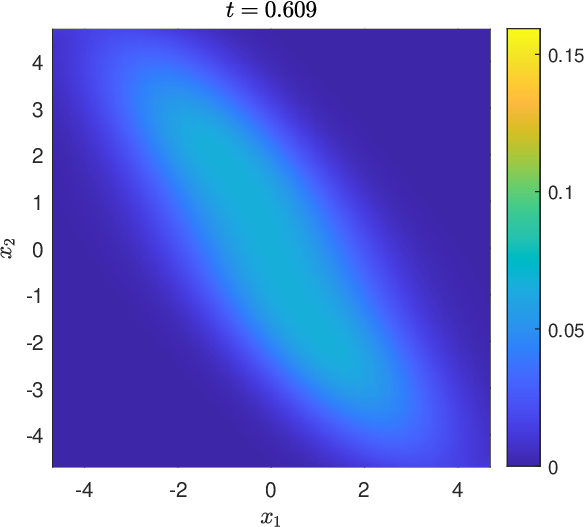}
    \label{fig4d}
}
\caption{
Numerical solution for the example in Section~\ref{sec:2d_local} with $\delta_n=n^{-1}$ after $160$ iterations.
Snapshots of the updated density $\bar{m}(x_1,x_2,t)$ at four time instants.
The initial bimodal distribution (\textbf{Top-Left}) merges toward the origin under the inward drift, and a fraction of agents exits as congestion near the center increases.
}
\label{fig4}
\end{figure}

Figure~\ref{fig4} shows snapshots of the updated density $\bar{m}$ obtained with $\delta_n=n^{-1}$ after $160$ iterations. Starting from the two peaks at $x_1=\pm 1.5$ (Figure~\ref{fig4a}), the inward drift concentrates mass toward the origin. Congestion then triggers exit, so the remaining density at $t=T$ is lower and concentrated near $x=(0,0)$ (Figure~\ref{fig4d}).

\begin{figure}
\centering
\subfigure[$m_{n}(x,T)$, $\delta_n\!=\!1$, odd $n$]{
    \centering
    \includegraphics[width=0.48\linewidth]{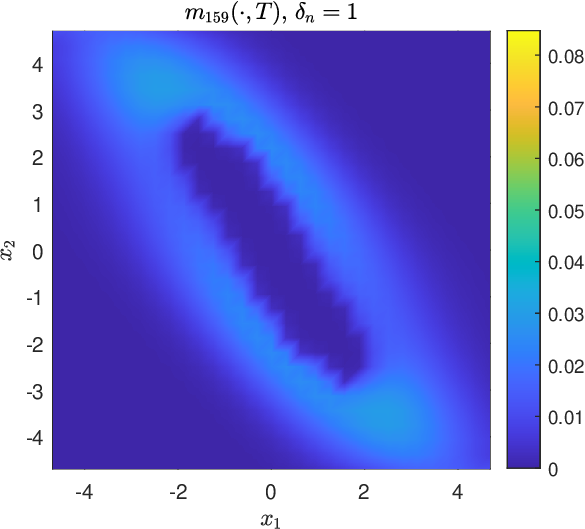}
    \label{fig5a}
}
\subfigure[$m_{n}(x,T)$, $\delta_n\!=\!1$, even $n$]{
    \centering
    \includegraphics[width=0.48\linewidth]{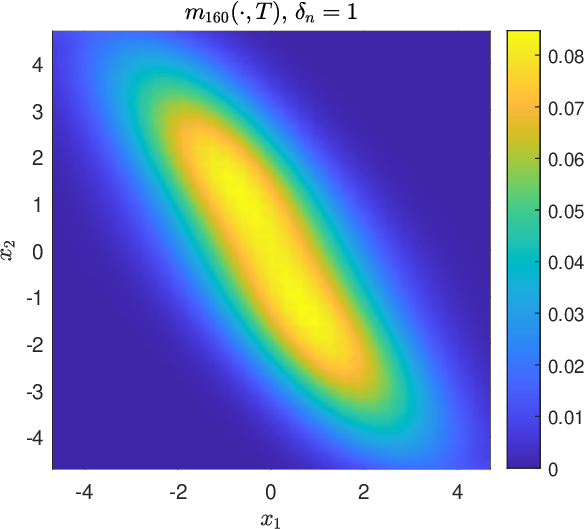}
    \label{fig5b}
}\\
\subfigure[$\bar{m}_{n}(x,T)$, $\delta_n\!=\!n^{-1}$, final iter.]{
    \centering
    \includegraphics[width=0.48\linewidth]{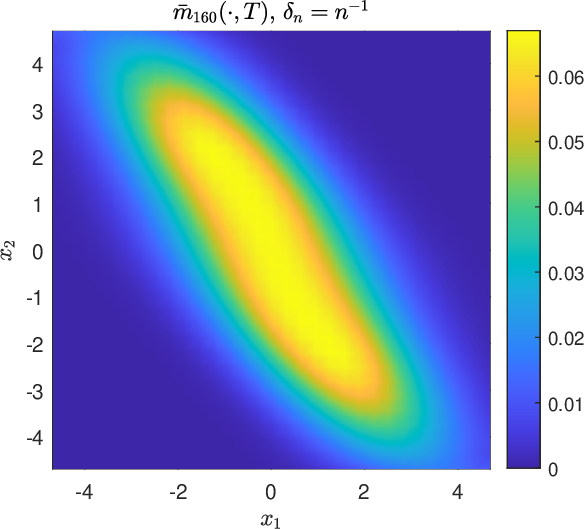}
    \label{fig5c}
}
\subfigure[Terminal mass of $\bar{m}_n$ vs.\ iteration]{
    \centering
    \includegraphics[width=0.48\linewidth]{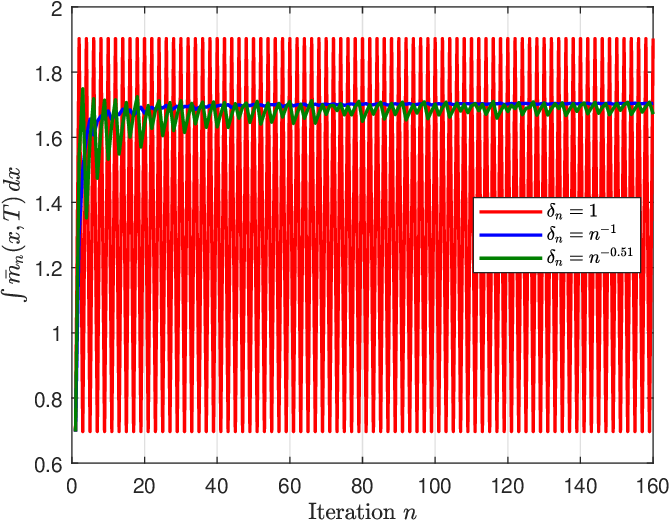}
    \label{fig5d}
}
\caption{
Oscillation vs.\ averaging for the example in Section~\ref{sec:2d_local}.
\textbf{Top-Left} and \textbf{Top-Right}: proposed terminal density $m_n(\cdot,T)$ at two consecutive iterations with $\delta_n=1$. The pure-strategy map is trapped in a cycle-$2$ orbit.
\textbf{Bottom-Left}: updated terminal density $\bar{m}_n(\cdot,T)$ after $160$ iterations with $\delta_n=n^{-1}$.
\textbf{Bottom-Right}: terminal mass $\int\bar{m}_n(x,T)\,dx$ versus iteration for $\delta_n=1$, $\delta_n=n^{-1}$, and $\delta_n=n^{-0.51}$.
}
\label{fig5}
\end{figure}

Figure~\ref{fig5} compares the pure-strategy map with fictitious-play averaging. With $\delta_n=1$, one has $\bar{m}_n=m_n$, and the terminal density alternates between two distinct profiles (Figures~\ref{fig5a} and~\ref{fig5b}); the corresponding terminal mass jumps between two values and does not settle (Figure~\ref{fig5d}). Thus the best-reply map of the pure-strategy system has no stable fixed point. With $\delta_n=n^{-1}$, which satisfies~\eqref{cond_delta}, the updated distribution $\bar{m}_n$ approaches a single profile (Figure~\ref{fig5c}) and its terminal mass stabilizes. The associated proposed distributions $m_n$ continue to vary from one iteration to the next: averaging a sequence of non-convergent pure strategies produces a mixed strategy equilibrium. The more aggressive admissible rate $\delta_n=n^{-0.51}$ also damps the cycle-$2$ oscillation of $\bar{m}_n$, but more slowly than $\delta_n=n^{-1}$.

\begin{figure}
\centering
\subfigure[$l$-2 error]{
    \centering
    \includegraphics[width=0.48\linewidth]{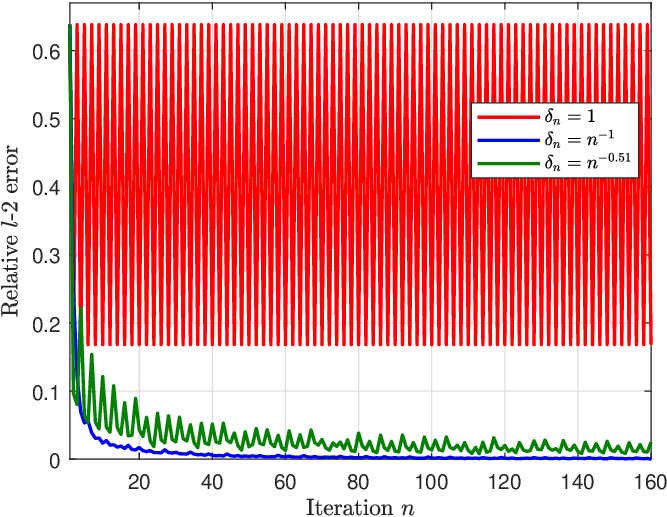}
    \label{fig6a}
}
\subfigure[Exploitability]{
    \centering
    \includegraphics[width=0.48\linewidth]{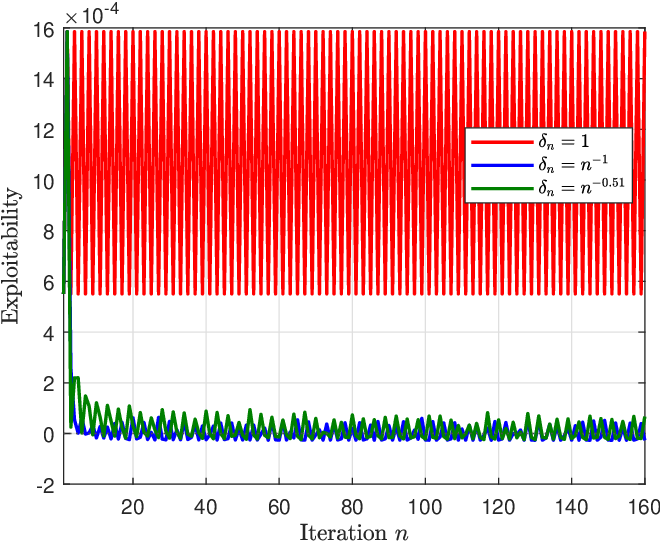}
    \label{fig6b}
}
\caption{
\centering
Relative $l$-2 error and exploitability versus iteration for the example in Section~\ref{sec:2d_local}.
\textbf{Left}: relative $l$-2 error of $\bar{m}_n$ with respect to the $\delta_n=n^{-1}$ result after $160$ iterations.
\textbf{Right}: exploitability.
With $\delta_n=1$ both quantities remain bounded away from zero and oscillate. With $\delta_n=n^{-1}$ and $\delta_n=n^{-0.51}$, which satisfy~\eqref{cond_delta}, the error of $\bar{m}_n$ decreases toward the baseline, while the exploitability becomes small; the decay is not monotone.
}
\label{fig6}
\end{figure}

For this problem the learning-rate condition in Definition~\ref{fictitiousplay} is essential. Figure~\ref{fig6} reports $160$ iterations for $\delta_n=1$, $\delta_n=n^{-1}$, and $\delta_n=n^{-0.51}$. The last two rates satisfy~\eqref{cond_delta}; $\delta_n=1$ does not. In the latter case the algorithm remains on the cycle-$2$ orbit of the pure-strategy map. In the former case $\bar{m}_n$ approaches a mixed strategy equilibrium, with $\delta_n=n^{-1}$ the more efficient of the two admissible rates tested here. As in Section~\ref{sec:nonlocal},~\eqref{cond_delta} is only sufficient: the point is that a rate violating it may fail when no pure strategy equilibrium exists, whereas a rate satisfying it yields a convergent fictitious-play sequence.

\section{Conclusion}

In conclusion, this paper proposes a novel generalized fictitious play algorithm for computing mixed strategy equilibria in OSMFGs. The key innovations include leveraging an iterative process of solving pure strategy systems to approximate mixed equilibria, as well as expanding the design flexibility for the learning rate parameter. Rigorous convergence results are provided, and finite difference schemes are constructed to efficiently solve the obstacle and Fokker-Planck equations during each iteration. 

Future work includes extensions to problems with common noise, where the equilibria consist of randomized stopping times that depend on the realized common noise path. The generalized fictitious play framework could also be applied to other competitive games involving optimal stopping decisions. Additionally, further analysis on quantifying the convergence rate and computational complexity could provide deeper theoretical insights. Overall, this paper introduces a novel algorithm and analysis to overcome the limitations of current OSMFG methods, opening the door for handling broader classes of large-scale dynamic games with optimal stopping.

\section*{Acknowledgement}
Jiajun Tong is partially supported by the Peking University Start-up Grant. Zhennan Zhou is supported by the National Key R\&D Program of China, Project Number 2021YFA1001200.
We thank Xu'an Dou and Jian-Guo Liu for helpful discussions.

\bibliographystyle{siam}
\bibliography{reference}
\end{document}